\documentclass[reqno]{amsart}
\usepackage{amssymb}
\usepackage{amsmath}

\makeatletter
\@addtoreset{equation}{section}
\makeatother

\renewcommand\thefigure{\thesection.\@arabic\c@figure}
\renewcommand\thetable{\thesection.\@arabic\c@table}
 
\newtheorem{theorem}{Theorem}[section]
\newtheorem{lemma}[theorem]{Lemma}
\newtheorem{proposition}[theorem]{Proposition}

\newtheorem{remark}[theorem]{Remark}
 
\newcommand{\mc}[1]{{\mathcal #1}}   
\newcommand{\mf}[1]{{\mathfrak #1}}
\newcommand{\mb}[1]{{\mathbf #1}}
\newcommand{\bb}[1]{{\mathbb #1}}

\newcommand{\id}{{1 \mskip -5mu {\rm I}}}
\renewcommand{\epsilon}{\varepsilon}
\renewcommand{\tilde}{\widetilde}

\newcommand{\<}{\langle}
\renewcommand{\>}{\rangle}

\begin{document}

\author{L.\ Bertini, A.\ De Sole, D.\ Gabrielli, G.\ Jona-Lasinio,
C.\ Landim} 

\address{Dipartimento di Matematica, Universit\`a di Roma La Sapienza, 
Piazzale Aldo Moro 2, 00185 Roma, Italy
\newline
e-mail:  \rm \texttt{bertini@mat.uniroma1.it}}

\address{Department of Mathematics, Harvard
  University, Cambridge, MA 02138, USA
\newline
e-mail:  \rm \texttt{desole@math.harvard.edu}}

\address{Dipartimento di Matematica, Universit\`a dell'Aquila,
  Coppito, 67100 L'Aquila, Italy
\newline
e-mail:  \rm \texttt{gabriell@univaq.it}}

\address{Dipartimento di Fisica and INFN, Universit\`a di Roma ``La
  Sapienza", Piazzale A. Moro 2, 00185 Roma, Italy 
\newline 
e-mail:  \rm \texttt{gianni.jona@roma1.infn.it}}

\address{\noindent IMPA, Estrada Dona Castorina 110,
CEP 22460 Rio de Janeiro, Brasil and CNRS UMR 6085,
Universit\'e de Rouen, 76128 Mont Saint Aignan, France.
\newline
e-mail:  \rm \texttt{landim@impa.br}}

\title[Current large deviations]
{Large deviations of the empirical current \\ ~ \\
 in interacting particle systems} 

\begin{abstract}
We study current fluctuations in lattice gases in the hydrodynamic
scaling limit.  More precisely, we prove a large
deviation principle for the empirical current in the symmetric simple
exclusion process with rate functional $I$.  We then estimate the
asymptotic probability of a fluctuation of the average current over a
large time interval and show that the corresponding rate function can
be obtained by solving a variational problem for the functional $I$.
For the symmetric simple exclusion process the minimizer is time
independent so that this variational problem can be reduced to a time
independent one. On the other hand, for other models the minimizer
is time dependent.  This phenomenon is naturally interpreted as
a dynamical phase transition.
\end{abstract}

                                           

\maketitle

\section{Introduction}
\label{sec0}     

In the last twenty years, interacting particle systems have become a
main subject of research in physics for the insight they provide on
the dynamical aspects of statistical physics.
On the mathematical side they provide a source of new
interesting problems in probability theory. Particularly relevant are
the results obtained on their hydrodynamical limits and the associated
large deviations since a microscopic derivation of
phenomelogical macroscopic laws can be rigorously established.  
More precisely, for symmetric conservative interacting particle
systems, it has been shown that the empirical density satisfies a
parabolic evolution equation.
The associated dynamical large deviations rate function measures the
asymptotic probability, as the number of particles diverges, of
fluctuations from the hydrodynamical evolution. 
As discussed in \cite{BDGJL2}, this rate function provides a new
approach to the analysis of stationary non equilibrium states.
These states describe a physical situation in which there is a
macroscopic flow through the system and the Gibbsian description 
is not applicable.  
Rigorous proofs of the dynamical large deviation principle have been
obtained for some equilibrium models, see e.g.\ \cite{dv,kov}, and for
the non equilibrium simple exclusion process \cite{bdgjl3}.

Beside the empirical density, a very important observable is the
current, which measures the flux of particle.  This quantity gives
informations that cannot be recovered from the density because from a
density trajectory we can determine the current trajectory only up to
a divergence free vector field. In \cite{BDGJL5,bl6} we have
introduced, at a heuristic level, the large deviation principle for
the empirical current.  In the present paper we prove it in the case
of the symmetric simple exclusion and illustrate some relevant
applications.

The simple exclusion process is a lattice gas with an exclusion
principle: a particle can move to a neighboring site only if this is
empty. The particle dynamics is given by a Markov process on the 
state space $\{0,1\}^{\bb T_N^d}$, where $\bb T_N^d = (\bb Z / N\bb
Z)^d$ is the discrete $d$-dimensional torus with $N^d$ points.
We denote by $\eta = \{\eta (x), x\in \bb T^d_N\}$ a configuration, so
that  $\eta(x)=1$ when the site $x$ is occupied, $\eta(x)=0$ otherwise. 
Let $\pi^N$ the empirical density of particles. The hydrodynamic limit
for this model is particularly simple: in the limit $N\to \infty$ 
the empirical density $\pi^N$ satisfies the heat equation.
To discuss the large deviation asymptotics we need to introduce the 
\emph{mobility} $\chi(\pi)$, which describes the response to 
an external field; for the symmetric simple exclusion we have 
$\chi(\pi)=\pi(1-\pi)$.
We introduce the integrated empirical current ${\mb W}_t^N$, 
which measures the total net flow of particles in the time interval $[0,t]$, 
associated to a trajectory $\eta$.  
We shall prove a large deviation principle which can be informally
written as follows. Fix a possible path $\mb W_t$, $t\in [0,T]$ of the
integrated empirical current, then   
\begin{equation}
\label{f1a}
\bb P^N \big( {\bf W}^N_t  \approx {\mb W}_t\,,\: t\in [0,T] \big) 
\sim \exp\big\{ - N^d \, I_{[0,T]}({\mb W} )\big\}
\end{equation}
where the rate functional is 
\begin{equation}
\label{Ica}
I_{[0,T]}( {\mb W} )\;=\; \frac 12 \int_0^T \!dt \,
\Big\langle \big[ {\dot {\mb W}_t} + \frac 12 \nabla \pi_t  \big], 
\frac{1}{\pi_t(1-\pi_t)}
\big[ {\dot {\mb W}_t}  + \frac 12 \nabla \pi_t \big] \Big\rangle
\end{equation}
In the above formula ${\dot {\mb W}_t}$ is the instantaneous current
at time $t$. Moreover $\pi_t$, which represents the associated
fluctuation of the empirical density, is obtained from ${\dot {\mb
W}}$ by solving the continuity equation $\partial_t\pi + \nabla \cdot
{\dot {\mb W}} =0$. Finally $\langle \cdot,\cdot\rangle$ denotes
integration with respect to the space variables.  Note that equation
\eqref{Ica} can be interpreted, in analogy to the classical Ohm's law,
as the total energy dissipated in the time interval $[0,T]$ by the
extra current ${\dot {\mb W}} + \frac 12 \nabla \pi$.  The large
deviation principle of the empirical density \cite{kov}, as we show,
can also be easily deduced from \eqref{f1a}--\eqref{Ica}.

Using \eqref{f1a}--\eqref{Ica} we then analyze the fluctuations properties
of the mean empirical current ${\mb W}^N_T/T = T^{-1} \int_0^T\!dt \:
\dot{\mb W}_t$ over a large time interval $[0,T]$.
This is the question addressed in \cite{bd} in one
space dimension by postulating an ``additivity principle'' which
relates the fluctuation of the current in the whole system to the
fluctuations in subsystems.  We show that the probability of observing
a given, divergence free, time averaged fluctuation $\mb J$ can be described
by a rate functional $\Phi(\mb J)$, i.e.\ as $N\to \infty$ and $T\to
\infty$ we have 
$$
\bb P^N \Big( \frac {{\mb W}^N_T}{T}  \approx {\mb J} \Big) 
\sim \exp\big\{ - N^d \, T  \, \Phi(\mb J ) \big\}
$$
The functional $\Phi$ is characterized
by a variational problem for the functional $I_{[0,T]}$
\begin{equation}
\label{limTa}
\Phi ({\bf J})
 = \lim_{T\to\infty} \; \inf_{{\mb W}} 
\frac 1T \; I_{[0,T]} ({\mb W})
\end{equation}
where the infimum is carried over all paths ${\mb W}_t$ 
such that $ {\mb W}_T = T {\bf J}$.

Let us denote by $U$ the functional obtained by restricting the
infimum in \eqref{limTa} to the paths $\mb W$ such that $\dot{\mb
W}_t$ is divergence free for any $t\in [0,T]$. The associated density
profile $\pi_t$ does not evolve. We get   
\begin{equation}
\label{Ua}
U({\mb J}) =  \inf_{\rho} \frac 12 
\Big\langle \big[ {\bf J} + \frac 12  \nabla \rho)  \big], 
\frac{1}{\rho(1-\rho)} \big[ {\bf J} + \frac 12 \nabla \rho \big] 
\Big\rangle
\end{equation} 
where the infimum is carried out over all the density profiles
$\rho$.  This is the functional introduced in \cite{bd} 
in one space dimension. 
For the symmetric simple exclusion process we prove 
that the additivity principle postulated in \cite{bd}
gives the correct answer, that is  $\Phi= U$.
On the other hand, while $\Phi$ is always convex the functional $U$ 
may be non convex. 
In general \cite{BDGJL5,bl6} we interpret the strict inequality
$\Phi<U$, as a dynamical phase transition. In such a case the
minimizers for \eqref{limTa} become in fact time dependent and the
invariance under time shifts is broken.   
In \cite{bl6} we have shown that, for the one-dimensional
Kipnis-Marchioro-Presutti (KMP) model \cite{bgl,kmp}, which is defined
by a harmonic chain with random exchange of energy between neighboring
oscillators, the following holds when it is considered with periodic
boundary condition.
The functional $U$ is given by $U(J)= (1/2) J^2/ \chi(m)= (1/2)
J^2/m^2 $, where $m$ is the (conserved) total energy. Moreover, for
$J$ large enough, $\Phi(J) < U(J)$.  This inequality is obtained by
constructing a suitable traveling wave current path whose cost is
less than $U(J)$.  In the present paper we give a more formal
presentation of these results.  We finally mention that the strict
inequality $\Phi < U$ also occurs for the weakly asymmetric exclusion
process for sufficiently large external field \cite{bd2}.

\section{Notation and results}
\label{sec1}

For $N\ge 1$, let $\bb T_N^d = (\bb Z / N\bb Z)^d$ be the discrete
$d$-dimensional torus with $N^d$ points. Consider the symmetric simple
exclusion process on $\bb T_N^d$. This is the Markov process on
the state space $\mc X_N:=\{0,1\}^{\bb T_N^d}$ whose generator 
$L_N$ is given by
\begin{eqnarray*}
(L_N f)(\eta) \;=\; \frac {N^2} 2 
\sum_{\substack{x,y\in\bb T_N^d \\ |x-y|=1}} 
\{ f(\sigma^{x,y} \eta)-f(\eta)\} \;.
\end{eqnarray*}
In this formula, $\eta = \{\eta(x),\, x\in\bb T_N^d\}\in \mc X_N $ is a
configuration of particles, so that $\eta (x)=0$, 
resp.\ $\eta(x) = 1$, if and only if site $x$ is
vacant, resp.\ occupied, for $\eta$, and $\sigma^{x,y}\eta$ is the
configuration obtained from $\eta$ by exchanging the occupation
variables $\eta(x)$, $\eta(y)$: 
$$
(\sigma^{x,y} \eta) (z)\; =\;
\left\{
\begin{array}{ll}
\eta (z)  & \hbox{if $z\neq x,y$}\; , \\
\eta (y)  & \hbox{if $ z=x$}\; , \\
\eta (x)  & \hbox{if $ z=y$}\; .
\end{array}
\right.
$$
Notice that we speeded up time by $N^2$. Denote by $\{\eta_t : t\ge
0\}$ the Markov process with generator $L_N$ and by $\bb P^N_{\mu^N}$
the probability measure on the Skorohod space $D(\bb R_+ , \mc X_N)$
induced by the Markov process $\eta_t$ and a probability measure
$\mu^N$ on $\mc X_N$, standing for the initial distribution. When
$\mu^N$ is a Dirac measure concentrated on a configuration $\eta^N$,
we denote $\bb P^N_{\mu^N}$ by $\bb P^N_{\eta^N}$.
An elementary computation shows that this process is reversible with
respect to any Bernoulli product measure on $\{0,1\}^{\bb T_N^d}$ with
parameter $m\in [0,1]$.

Denote by $\mc M (\bb T^d)$ the space of finite signed measures on
$\bb T^d$, the $d$ dimensional torus of side 1, 
endowed with the weak topology and by $\mc F = \mc
F_{+,1}(\bb T^d)$ the set of positive, measurable functions bounded by
$1$ endowed with the same weak topology. For a finite signed measure
$m$ we let $\< m, F\>$ be the
integral of a continuous function $F: \bb T^d \to\bb R$ with respect to
$m$. Likewise  for a profile $\pi\in \mc F$ and $F\in C(\bb T^d)$  
we denote by $\< \pi, F\>$ the integral of $\pi F$.

For $N\ge 1$ and a configuration $\eta\in \mc X_N$, denote by
$\pi^N=\pi^N(\eta^N)$ the empirical density of particles. It is defined as
$$
\pi^N (u)  \;=\; \sum_{x\in\bb T_N^d} \mb 1\{ u N \in B(x) \}
\, \eta (x) \;,
$$
where, for $x = (x_1, \dots, x_d)\in \bb T_N^d$, $B(x)$ is the cube 
$[x_1,x_1+1) \times \cdots \times [x_d, x_d +1)$. Set $\pi^N_t =
\pi^N(\eta_t)$ and notice that $\pi^N_t$ belongs to $\mc F$ for each
$t\ge 0$.

For $t\ge 0$ and two neighboring sites $x,y\in \bb T^d_N$, 
denote by $J^{x,y}_t$ the total number of particles that
jumped from $x$ to $y$ in the macroscopic time interval $[0,t]$.  
Let $\{e_k: 1\le k\le d\}$ be the canonical basis of $\bb R^d$; 
the difference $W^{x,x+e_j}_t = J^{x,x+e_j}_t - J^{x+e_j,x}_t$ represents
the net flow of particles across the bond $\{x,x+e_j\}$ in the time interval
$[0,t]$.

For $t\ge 0$, we define the \emph{empirical integrated current} 
$\mb W^N_t =(W^{N}_{1,t}, \dots, W^{N}_{d,t}) \in \mc M_d = \{\mc M(\bb
T^d)\}^d$ as the vector-valued finite signed measure on $\bb T^d$
induced by the net flow of particles in the time interval $[0,t]$:
$$
W^{N}_{j,t} \;=\; N^{-(d+1)} \sum_{x\in\bb T_N^d} 
W^{x,x+e_j}_t \delta_{x/N}\;,
\qquad j=1,\dots,d
$$ 
where $\delta_u$ stands for the Dirac measure concentrated on $u$.
Notice the extra factor $N^{-1}$ in the normalizing constant which
corresponds to the diffusive rescaling of time. In particular, for a
continuous vector field $\mb F = (F_1, \dots, F_d) \in C (\bb T^d;\bb
R^d)$ the integral of $\mb F$ with respect to $\mb W^N_t$, also
denoted by $\<\mb W^N_t, \mb F\> $, is given by
\begin{equation}
\label{empcurr}
\< \mb W^N_t , \mb F\> \;=\; N^{-(d+1)} \sum_{j=1}^d \sum_{x\in\bb T_N^d} 
F_j (x/N) \, W^{x,x+e_j}_t \;.
\end{equation}

The purpose of this article is to prove a large deviations principle
for the empirical integrated current $\mb W^N_t$ and discuss the
asymptotic behavior as $t\to \infty$. We start with the law of large
numbers. Fix a profile $\lambda \in \mc F$ and let $\{\mu^N : N\ge
1\}$ a sequence of measures on $\mc X_N$ associated to $\lambda$ in
the sense that the empirical density converges to $\lambda$ in
probability with respect to $\mu^N$.  Namely, for each $F\in C(\bb
T^d)$ and $\delta>0$, we have
\begin{equation}
\label{f01}
\lim_{N\to\infty} \mu^N \Big\{ \, \Big\vert \<\pi^N, F\> 
\;-\; \int_{\bb T^d} \lambda (u) F(u) \, du \Big\vert > \delta \Big\}
\;=\;0\;. 
\end{equation}
It is well known, see e.g.\ \cite{kl}, that in such a case the empirical
density $\pi^N_t$ converges in probability to $\rho=\rho(t,u)$ which
solves the heat equation
\begin{equation}
\label{f02}
\left\{
\begin{array}{l}
\partial_t \rho \;=\; (1/2) \Delta \rho \;, \\
\rho(0, \cdot) = \lambda (\cdot) \;,
\end{array}
\right.
\end{equation}
where $\Delta = \nabla \cdot \nabla$ stands for the Laplacian and
$\nabla$ for the gradient.  We claim that the empirical current
converges to the time integral of $- (1/2) \nabla \rho$. This is the
content of the next result which is proven in Subsection~\ref{sec2} in
a more general context.

\begin{proposition}
\label{s1}
Fix a profile $\lambda\in \mc F$ 
and consider a sequence of probability measures $\mu^N$ 
associated to $\lambda$ in the sense of \eqref{f01}. 
Let $\rho$ be the solution of the heat equation
\eqref{f02}. Then, for each $T>0$, $\delta>0$ and $\mb F \in C(\bb
T^d;\bb R^d)$,
\begin{equation*}
\lim_{N\to\infty} \bb P^N_{\mu^N} \Big[\, \Big\vert \<\mb W^N_T, \mb F\> 
\;+\; (1/2) \int_0^T dt\, \int_{\bb T^d} \mb F(u) \cdot 
\nabla \rho(t,u) \, du \Big\vert > \delta \Big ] \;=\; 0\; .
\end{equation*}
\end{proposition}

We turn now to the large deviations of the pair $(\mb W^N, \pi^N)$.  Fix
$T>0$ and denote by $D([0,T], \mc M_d \times \mc F)$ the space of
cad-lag trajectories with values in $\mc M_d \times \mc F$ endowed
with the Skorohod topology.  
Fix a profile $\gamma\in C^2(\bb T^d)$ bounded away from $0$ 
and $1$: there exists $\delta >0$ such that $\delta \le \gamma \le
1-\delta$.
To focus on the dynamical fluctuations, 
we assume that the process starts from a deterministic initial
condition $\eta^N$ which is associated to $\gamma$ in the sense that
$\pi^N(\eta^N) \to \gamma$ in $\mc F$.

Let $\mf A_\gamma$ be the set of trajectories
$(\mb W, \pi)$ in $D([0,T], \mc M_d \times \mc F)$ such that
for any $t\in[0,T]$ and any $F \in C^1(\bb T^d)$ 
\begin{equation}
\label{f09}
\langle \pi_t , F\rangle - \langle \gamma , F\rangle = 
\langle \mb W_t , \nabla F \rangle 
\,,\qquad \mb W_0 = 0
\end{equation}
Note that $\mf A_\gamma$ is a \emph{closed and convex} subset of 
$C([0,T], \mc M_d \times \mc F)$.
Equation \eqref{f09} is the weak formulation of the continuity equation 
$\dot{\pi}_t + \nabla \cdot \dot{\mb W}_t = 0$, 
$\pi_0=\gamma$, ${\mb W}_0 = 0$. 
For each $\mb F \in
C^{1,1}([0,T]\times \bb T^d; \bb R^d)$, 
define the \emph{convex and lower semi-continuous}
functional $J_{\mb F}$ as follows. If $(\mb W,\pi)\in C([0,T];\mc
M_d\times\mc F)$ we set 
\begin{eqnarray}
\label{djf}
J_{\mb F}(\mb W, \pi) &=& \< \mb W_T, \mb F_T\> \; 
-\; \int_0^T dt \, \< \mb W_t, \partial_t \mb F_t\> 
\nonumber\\ 
&-& \frac 12 \int_0^T dt \, \< \pi_t, \nabla \cdot \mb F_t\> 
\;-\; \frac 12 \int_0^T dt \, \< \chi(\pi_t), | \mb F_t| ^2 \>\; ,
\end{eqnarray}
where, here and below, $\chi(a) = a (1-a)$ is the mobility. We set 
$J_{\mb F}(\mb W, \pi) =+\infty$ if  $(\mb W,\pi)\not\in C([0,T];\mc
M_d\times\mc F)$. Since $ C([0,T];\mc M_d\times\mc F)$ is a closed
subset of  $D([0,T];\mc M_d\times\mc F)$, the argument in 
\cite[\S 10.1]{kl} proves the convexity and lower semi-continuity of
$J_{\mb F}$.

Let finally
\begin{equation}
\label{f10}
J(\mb W, \pi) = \sup_{\mb F \in C^{1,1}} J_{\mb F}(\mb W, \pi)
\;, \quad I(\mb W, \pi) \; =\;
\left\{
\begin{array}{ll}
J (\mb W, \pi) & \text{if $(\mb W, \pi) \in \mf A_\gamma$}\;, \\
+ \infty & \text{otherwise}\;.
\end{array}
\right.
\end{equation}

Notice that the functional $J$ is convex and lower semi-continuous,
properties which are inherited by $I$ because $\mf A_\gamma$ is a
closed and convex subset of $C([0,T], \mc M_d \times \mc F)$.  
Notice furthermore that the continuity equation
\eqref{f09} determines the trajectory $\pi$ as a function of $\mb W$
and the initial condition $\gamma$.  Therefore, the rate function $I
(\mb W, \pi)$ can be thought as a function of $\mb W$ and the initial
density profile.
In Subsection~\ref{sec4} we derive a more explicit formula for the rate
function \eqref{f10}. If $I(\mb W,\pi) < \infty$ we have
\begin{equation*}
I(\mb W, \pi) \; =\; \frac 12 \int_0^T dt \, \< \chi(\pi_t), | \mb
F_t| ^2 \>\; , 
\end{equation*}
where the vector-valued function $\mb F_t$ is the solution of
\begin{equation*}
\partial_t \mb W_t \;+\; (1/2) \nabla \pi_t \;=\; \chi(\pi_t) \mb F_t\;.
\end{equation*}
Thus, formally, 
\begin{equation*}
I(\mb W, \pi) \; =\; \frac 12 \int_0^T dt \, \Big\< \frac 1{\chi(\pi_t)}
\, \big | \dot{ \mb W}_t -  \dot{ \mb W}_t(\pi) \big | ^2 \Big\>\; , 
\end{equation*}
where $\dot{ \mb W}_t = \partial_t \mb W_t$ is the instantaneous
current at time $t$ for the path $(\mb W, \pi)$ and 
$\dot{ \mb W}_t(\pi) = - (1/2) \nabla \pi_t$ is the typical
instantaneous current at 
time $t$ associated to the density profile $\pi_t$. Recall in fact that
the hydrodynamic equation \eqref{f02} is in our case the heat equation. 

\begin{theorem}
\label{s5}
Consider a profile $\gamma\in C^2(\bb T^d)$ bounded away from $0$ and
$1$ and a sequence $\{\eta^N : N\ge 1\}$ such that
$\pi^N(\eta^N)\to\gamma$ in $\mc F$.  Then, for each closed set $F$
and each open set $G$ of $D([0,T], \mc M_d \times \mc F)$, we have
\begin{eqnarray*}
\!\!\!\!\!\!\!\!\!\!\!\!\! 
&&
\varlimsup_{N\to\infty} \frac 1{N^d} \log \bb P_{\eta^N}^N
\Big[ (\mb W^N, \pi^N) \in F \Big] \;\le\; - \inf_{(\mb W, \pi) \in F}
I (\mb W, \pi) \;, \\
\!\!\!\!\!\!\!\!\!\!\!\!\! && \qquad
\varliminf_{N\to\infty} \frac 1{N^d} \log \bb P_{\eta^N}^N
\Big[ (\mb W^N, \pi^N) \in G \Big] \;\ge\; - \inf_{(\mb W, \pi) \in G}
I (\mb W, \pi) \;.
\end{eqnarray*}
\end{theorem}

We next state a ``large deviation principle'' for the mean empirical
current $\mb W^N_T/ T$ in the interval $[0,T]$ as we let \emph{first}
$N\to\infty$ and \emph{then} $T\to\infty$.

Let us denote by $\mc B \subset \mc M_d$ the set of divergence free
measures, i.e.\  
\begin{equation}
\label{dB}
\mc B := \Big\{ {\mb J}\in \mc M_d \, : \, 
\langle {\mb J}, \nabla f \rangle = 0 
\; \text{ for any } f \in C^1(\bb T^d)  
\Big\}
\end{equation}
which is a closed subspace of $\mc M_d$.
Given $m\in (0,1)$, we introduce the set of profile with mass $m$,
i.e.\ we set $\mc F_m:= \big\{ \rho \in \mc F \,:\,
\int_{\bb T^d}\!du \: \rho(u) = m \big\}$. We finally define $U_m: \mc M_d \to
[0,+\infty]$ as 
\begin{equation}
\label{dU}
U_m({\mb J}) := 
\inf_{\substack{\rho\in \mc F_m\cap C^2(\bb T^d)\\ 0< \rho< 1 }} 
\; \frac 12 
\Big\langle \big[ {\bf j} + \frac 12 \nabla \rho \big], \frac 1{ \rho
    (1-\rho)} \big[ {\bf j} + \frac 12 \nabla \rho \big]
\Big\rangle   
\end{equation}
if ${\mb J} \in \mc B \;$, ${\mb J}(du) = {\mb j}\, du$, and  
$U_m({\mb J})= +\infty$  otherwise.  
In Section~\ref{s:p=u} we show that $U_m$ is a lower semi-continuous
convex functional. 

In Section~\ref{s:p=u} we also prove that in the one dimensional case, where
$J\in \mc B$, $J(du) =j du$, implies that $j$ is constant $du$-a.e., we simply
have 
\begin{equation}
\label{dud=1}
U_m(J) := 
\begin{cases}
{\displaystyle 
 \frac 12  \, \frac {j^2}{m(1-m)}
}
& \text{if $J=j du$  for some $j$ $du$-a.e. constant} \\ 
+\infty & \text{otherwise}
\end{cases}
\end{equation}

\begin{theorem}
\label{ldpnt}
Let $m\in (0,1)$, $\gamma\in C^2 (\bb T^d) \cap \mc F_m$ 
bounded away from $0$ and $1$, and 
$\eta^N\in \mc X_N$  a sequence such that $\pi^N(\eta^N)\to
\gamma$ in $\mc F$.
Then, for each closed set $C$ and each open set $G$ of $\mc M_d$,
we have
\begin{eqnarray*}
&&
\varlimsup_{T\to\infty} \varlimsup_{N\to\infty} 
\frac 1{T N^d} \log 
\bb P_{\eta^N}^N
\Big[ \frac 1T {\bf W}^N_T \in C \Big] \;\le\; 
- \inf_{{\bf J}\in C} U_m ({\bf J}) \;, 
\\ 
&&\qquad
\varliminf_{T\to\infty} \varliminf_{N\to\infty} 
\frac 1{T N^d} \log \bb P_{\eta^N}^N
\Big[ \frac 1T {\bf W}^N_T \in G \Big] \;\ge\; 
- \inf_{{\bf J}\in G} U_m({\bf J}) \;. 
\end{eqnarray*}
\end{theorem}

\section{Large deviation for the empirical current on a fixed time interval}
\label{s:dd}

In this Section we prove Theorem~\ref{s5}.
The proof is similar to the one of the large deviations principle for
the empirical density, see \cite{kov} or \cite[Chapter~10]{kl}. 
We therefore present only the main modifications.

\subsection{Weakly asymmetric exclusion processes}
\label{sec2}

In this subsection we prove the law of large numbers for the empirical
current of weakly asymmetric exclusion processes. Proposition~\ref{s1}
follows as a particular case.

Fix $T>0$ and a time dependent vector-valued function $\mb F =(F_1,
\dots, F_d) \in C^{1,1}( [0,T]\times \bb T^d ;\bb R^d)$. Denote by
$L_{\mb F,N}$ the time-dependent generator on $\mc X_N$ given by
\begin{eqnarray*}
(L_{\mb F,N} f)(\eta) \;=\; \frac {N^2} 2 
\sum_{j=1}^d \sum_{x\in\bb T_N^d} 
c^{\mb F}_{x,x+e_j} (t,\eta) 
\, \{ f(\sigma^{x,x+e_j} \eta)-f(\eta)\} \;,
\end{eqnarray*}
where the rate $c^{\mb F}_{x,x+e_j} (t,\eta)$ is given by
$$
\eta(x) [ 1- \eta(x+e_j)] \, e^{N^{-1} F_j(t, x/N)} 
\;+\; \eta(x+e_j) [ 1- \eta(x)] \, e^{- N^{-1} F_j(t, x/N)} \;.
$$
Hence, for $N$ large, instead of jumping from $x$ to $x+e_j$ (resp. $x+e_j$ to
$x$) with rate $1/2$, at a macroscopic time $t$ particles jump with
rate $(1/2) \{ 1 + N^{-1} F_j(t,x/N) \}$ (resp. $(1/2) \{ 1 - N^{-1}
F_j(t,x/N) \}$) and a small drift appears due to the external field
$\mb F$. For a probability measure $\mu^N$ on $\mc X_N$, denote by
$\bb P^N_{\mb F, \mu^N}$ the measure on the path space $D(\bb R_+ ,
\mc X_N)$ induced by the Markov process $\eta_t$ with generator
$L_{\mb F,N}$ and initial distribution $\mu^N$.

Let $\rho^{\mb F,\lambda}$ be the unique weak solution of the parabolic
equation
\begin{equation}
\label{f03}
\left\{
\begin{array}{l}
\partial_t \rho \;=\; (1/2) \Delta \rho - \nabla \cdot \{ \chi(\rho)
\mb F\} \;, \\
\rho(0, \cdot) = \lambda (\cdot) \;.
\end{array}
\right.
\end{equation}
Write the previous differential equation as
\begin{equation*}
\dot {\rho}_t \;+\; \nabla \cdot \dot{\mb W}_{\mb F}(\rho_t) \;=\; 0\; ,
\end{equation*}
where $\dot{\mb W}_{\mb F}(\rho)$ is the instantaneous current
associated to the profile $\rho$ and is given by
\begin{equation*}
\dot{\mb W}_{\mb F}(\rho) \;=\; - (1/2) \nabla \rho 
\;+\; \chi(\rho) \mb F\;.
\end{equation*}
The main result of this section states that $\mb W^N_t$ converges in
probability to the time integral of $\dot{\mb W}_{\mb F}(\rho)$:

\begin{lemma}
\label{s2}
Fix a profile $\lambda : \bb T^d\to[0,1]$ and consider a sequence of
probability measures $\{\mu^N : N\ge 1\}$ on $\mc X_N$ associated to
$\lambda$ in the sense of \eqref{f01}.  For each $t>0$, $\delta>0$,
$G\in C(\bb T^d)$ and $\mb H \in C^1( \bb T^d;\bb R^d)$,
\begin{eqnarray*}
\!\!\!\!\!\!\!\!\!\!\!\!\! &&
\lim_{N\to\infty} \bb P^N_{\mb F, \mu^N} \Big[\, \Big\vert 
\<\pi_t^N, G\> \;-\; \big \< \rho^{\mb F, \lambda}_t, G \big \>
\Big\vert > \delta \Big ] \;=\; 0\; , \\
\!\!\!\!\!\!\!\!\!\!\!\!\! && \quad
\lim_{N\to\infty} \bb P^N_{\mb F, \mu^N} \Big[\, \Big\vert 
\<\mb W_t^N, \mb H\> \;-\; \int_0^t ds\, 
\big \< \dot{\mb W}_{\mb F}(\rho^{\mb F, \lambda}_s), \mb H \big \>
\Big\vert > \delta \Big ] \;=\; 0\; ,
\end{eqnarray*}
where $\< \dot{\mb W}_{\mb F}(\rho^{\mb F, \lambda}_s), \mb H \>$ stand for
\begin{equation*}
(1/2) \<\rho^{\mb F, \lambda}_s , \nabla \cdot \mb H \> 
\;+\; \<\chi(\rho^{\mb F, \lambda}_s) , \mb F_s
\cdot \mb H \> \; .
\end{equation*}
\end{lemma}

\begin{proof}
The law of large numbers for the empirical density follows from the
usual entropy method, see e.g. \cite[Chapter~6]{kl}. 
For each $t\ge 0$ the empirical density $\pi^N_t$ converges in
probability to $\rho^{\mb F, \lambda} (t,\cdot)$.
  
To derive the hydrodynamic equation for the current, fix $t>0$ and a
smooth vector field $\mb H : \bb T^d \to \bb R^d$. Let $\widetilde
{\mb W}^{N,\mb H}_t$ be the martingale defined by
\begin{eqnarray*}
\widetilde {\mb W}^{N,\mb H}_t & \!\!\! = \!\!\! & 
\< \mb W^N_t , \mb H\> - 
\int_0^t ds \frac {N^2} {2N^{d+1}} \sum_{j, x} H_j (x/N)\, 
\eta_s(x) [ 1-\eta_s(x+e_j)] \, e^{F_j(s,x/N)/N}
\\
& \!\!\! + \!\!\! & \int_0^t ds \frac {N^2}{2N^{d+1}} \sum_{j, x} 
H_j (x/N) \, \eta_s(x+e_j) [ 1-\eta_s(x)] \,
e^{- F_j(s,x/N)/N}\;.
\end{eqnarray*}
An elementary computation shows that the quadratic variation of this
martingale vanishes in $L^1 (\bb P^N_{\mb F, \mu^N})$ as
$N\uparrow\infty$. On the other hand, after a Taylor expansion and
few summations by parts, the time integral can be rewritten as
\begin{eqnarray*}
\!\!\!\!\!\!\!\!\!\!\!\!\! &&
\int_0^t ds\, \big\<\pi^N_s, (1/2) \nabla \cdot \mb H + \mb F_s
\cdot \mb H \big \> \; +\; O_{\mb F, \mb H} (N^{-1}) \\
\!\!\!\!\!\!\!\!\!\!\!\!\! && \quad
\;-\; \int_0^t ds\,
\frac 1{N^d} \sum_{j=1}^d \sum_{x\in\bb T_N^d}
H_j(x/N) F_j(s,x/N) \eta_s(x) \eta_s(x+e_j) \;,
\end{eqnarray*}
where $O_{\mb F, \mb H} (N^{-1})$ is an expression whose absolute
value is  bounded by $C N^{-1}$ for some 
constant depending only on $\mb F$ and $\mb H$.
By the two block estimates and the law of large numbers for the
empirical density, as $N\uparrow\infty$, the previous expression
converges in $\bb P^N_{\mb F, \mu^N}$-probability to
\begin{equation*}
(1/2)  \int_0^t ds\, \big\<\rho^{\mb F, \lambda}_s , 
\nabla \cdot \mb H \big\> 
\;+\; \int_0^t ds\, \big\<\chi(\rho^{\mb F, \lambda}_s) , \mb F_s
\cdot \mb H \big \> \; .
\end{equation*}
Since the martingale $\widetilde {\mb W}^{N,\mb H}_t$ vanishes in
$L^2$ as $N\uparrow\infty$, the lemma is proven.
\end{proof}

The same result holds for the generator $\tilde L_{\mb F,N}$ defined
by 
\begin{equation*}
(\tilde L_{\mb F,N} f)(\eta) \;=\; \frac {N^2} 2
\sum_{\substack{x,y\in\bb T_N^d \\ |x-y|=1}}
\eta(x) [ 1- \eta(y)] \, e^{N^{-1} \mb F(t, x/N) \cdot (y-x)}
\, \{ f(\sigma^{x,y} \eta)-f(\eta)\} \;.
\end{equation*}
However, the computations of the exponential martingales in next
section are slightly more complicated if we use this expression
instead of $L_{\mb F,N}$.

\subsection{Large deviations upper bound}
\label{sec3}

We first remark that Lemma~\ref{s6} below implies that 
the probability of the event  $(\mb W^N,\pi^N)\not \in C([0,T];\mc
M_d\times\mc F)$ is super-exponentially small as
$N\to\infty$. Recalling the definition \eqref{f10} of the rate
function $I$, it is therefore enough to prove the upper bound for
closed subsets of $C([0,T];\mc M_d\times\mc F)$. We shall first prove
it for compacts and then show the exponential tightness.

We start by recalling the super-exponential estimate of
\cite{dv,kov}. 
For a positive integer $\ell$ and $x$ in $\bb Z^d$, denote by $\eta^\ell
(x)$ the empirical density of particles on a box of size $2\ell +1$
centered at $x$:
$$
\eta^\ell(x)\;=\; \frac 1{(2\ell +1)^d}
\sum_{|y-x | \le \ell} \eta(y)\;.
$$
Moreover, for $1\le j\le d$, $\epsilon>0$, let
$$
V_{j, N,\epsilon} (\eta) \;=\; \frac 1{N^d} \sum_{x\in\bb T_N^d}
\Big | \frac 1{(2 \epsilon N +1)^d} \sum_{|y-x|\le \epsilon N} 
\eta(y) \eta(y+e_j) -
\big[\eta^{N\epsilon}(x)\big]^2 \Big|\; .
$$
     
\begin{theorem}
\label{s3}
For each $1\le j\le d$, $T>0$, each sequence of measures $\{ \mu^N :
N\ge 1\}$ and each $\delta>0$,
$$
\limsup_{\epsilon\to 0}\,\limsup_{N\to\infty} \frac 1{N^d}  
\log \bb P^N_{\mu^N} \Big[\, \Big | \int_0^{T} V_{j, N,\epsilon}(\eta_t) 
\, dt \Big | \;>\; \delta \Big] \;=\; -\infty \; .
$$
\end{theorem}

Fix a vector-valued function $\mb F : [0,T]\times \bb T^d \to \bb R^d$
in $C^{1,1}$ and denote by $d \bb P^N_{\mb F, \mu^N} /d \bb
P^N_{\mu^N} (T)$ the Radon-Nikodym derivative of $\bb P^N_{\mb F,
  \mu^N}$ with respect to $\bb P^N_{\mu^N}$ restricted to the time
interval $[0,T]$. A long but elementary computation gives that
\begin{eqnarray*}
\!\!\!\!\!\!\!\!\!\!\!\! &&
\frac 1{N^d} \log \frac {d \bb P^N_{\mb F, \mu^N} }{d \bb P^N_{\mu^N}}
(T) \;=\; \frac 1{N^{d+1}} \sum_{j=1}^d\sum_{x\in\bb T_N^d} \int_0^T
F_j(t,x/N) dW_t^{x,x+e_j} \\
\!\!\!\!\!\!\!\!\!\!\!\! &&
- \frac 12 \int_0^T dt \< \pi^N_t, \nabla \cdot \mb F_t\>
- \frac 14 \int_0^T dt \frac 1{N^d} \sum_{j=1}^d\sum_{x\in\bb T_N^d}
F_j(t,x/N)^2 \tau_x h_j(\eta_t) \;+\; O_{\mb F}(N^{-1})\; ,
\end{eqnarray*}
where $h_j (\eta) = \eta(0) + \eta(e_j) - 2\eta(0) \eta(e_j)$ and
$\tau_x$ denotes the translation of $x$.
In particular, on the set 
\begin{equation*}
\sum_{j=1}^d \Big | \int_0^{T} V_{j, N,\epsilon}(\eta_t) 
\, dt \Big | \;\le\; \delta\;,
\end{equation*}
integrating by parts the first term on the right hand side of the
penultimate formula, we obtain that $N^{-d} \log \{d \bb P^N_{\mb F,
  \mu^N} /d \bb P^N_{\mu^N}\} (T)$ is bounded below by
\begin{equation*}
J_{\mb F, \epsilon , \delta}(\mb W^N, \pi^N) \;-\; 
C(\mb F) \{ \epsilon + \delta \}
\end{equation*}
for every $\epsilon >0$. Here,
\begin{eqnarray*}
J_{\mb F, \epsilon , \delta}(\mb W, \pi) &=&
\< \mb W_T, \mb F_T\> \; -\; \int_0^T dt \, \< \mb W_t, \partial_t \mb
F_t\> \\ 
&-& \frac 12 \int_0^T dt \< \pi_t, \nabla \cdot \mb F_t\> 
\;-\; \frac 12 \int_0^T dt \, \< \chi(\pi^{\epsilon}_t),
| \mb F_t| ^2 \>\; ,
\end{eqnarray*}
$C(\mb F)$ is a finite constant depending only on $\mb F$ and
$\pi^{\epsilon}_t$ is the function defined by $\pi^{\epsilon}_t(u) =
(2\epsilon)^{-d} \int_{[u-\epsilon \mb 1, u+\epsilon \mb 1]} dv
\pi_t(v)$, where $[u-\epsilon \mb 1, u+\epsilon \mb 1]$ is the
hyper-cube $[u_1-\epsilon, u_1+\epsilon] \times \cdots \times
[u_d-\epsilon, u_d+\epsilon]$.

For each $\epsilon >0$, $\delta>0$ and $\mb F$ of class $C^{1,1}$, 
the functional  $J_{\mb F, \epsilon , \delta}$ is continuous  in $C([0,T];\mc
M_d\times\mc F)$. By repeating the
arguments presented in\cite[\S~10.4 ]{kl}, we then obtain that for each
compact set $K$ of $C([0,T], \mc M_d \times \mc F)$,
\begin{equation*}
\limsup_{N\to\infty} \frac 1{N^d} \log \bb P^N_{\eta^N} \Big[\, (\mb
W^N, \pi^N) \in K \Big] \;\le\; -\inf_{(\mb W, \pi) \in K} J (\mb W, \pi)
\; ,
\end{equation*}
where $J$ is defined in \eqref{f10}.

\medskip
To extend the upper bound from compact to closed sets, we next prove
the exponential tightness of the sequence $(\mb W^N,\pi^N)$.  As
stated before, the following Lemma also implies that the probability
of the event $(\mb W^N,\pi^N)\not \in C([0,T];\mc M_d\times\mc F)$ is
super-exponentially small.

\begin{lemma}
\label{s6}
Fix a sequence of measures $\{ \mu^N : N\ge 1\}$, a continuous
function $G: \bb T^d\to\bb R$, a vector-valued function $\mb H: \bb
T^d\to\bb R^d$ in $C^1$ and $\epsilon >0$. Then,
\begin{eqnarray*}
\!\!\!\!\!\!\!\!\!\!\! &&
\lim_{\delta\to 0} \limsup_{N\to\infty} \frac{1}{N^d} \log \bb
P^N_{\mu^N} \Big[ \sup_{|t-s|\le \delta} \big\vert <\pi_t^N, G> -
<\pi_s^N, G>\big\vert\; >\; \epsilon \Big ]\; =\; -\infty\; , \\
\!\!\!\!\!\!\!\!\!\!\! && \quad
\lim_{\delta\to 0} \limsup_{N\to\infty} \frac{1}{N^d} \log \bb
P^N_{\mu^N} \Big[ \sup_{|t-s|\le \delta} \big\vert <\mb W^N_t, \mb H> -
<\mb W^N_s, \mb H>\big\vert\; >\; \epsilon \Big ]\; =\; -\infty\; .
\end{eqnarray*}
\end{lemma}

\begin{proof}
The proof of the first estimate is similar to the one of the
exponential tightness of the empirical measure presented in
\cite[\S~10.4]{kl}. In our context the initial configuration is not
however the invariant measure.  The necessary modifications are 
worked out below in the proof of the second statement of the lemma.

To prove the second estimate, first observe that by a triangular
inequality and since
\begin{equation}
\label{f13}
\limsup_{N\to\infty} N^{-d} \log
\{a_N + b_N\} \le \max\Big \{ \limsup_{N\to\infty} N^{-d} \log a_N ,
\limsup_{N\to\infty} N^{-d} \log b_N \Big \}\;,
\end{equation}
it is enough to estimate
\begin{equation}
\label{f12}
\max_{0\le k\le T\delta^{-1}} \limsup_{N\to\infty} \frac{1}{N^d} \log \bb
P^N_{\mu^N} \Big[ \sup_{t_k\le t\le t_{k+1}} \big\vert <\mb W_t^N, \mb H> -
<\mb W_{t_k}^N, \mb H>\big\vert\; >\; \epsilon/3 \Big ]\;,
\end{equation}
where $t_k = k\delta$. By \eqref{f13}, we may also disregard the
absolute value in the previous expression provided we estimate the same
term with $- \mb H$ in place of $\mb H$. Fix $a>0$ and denote by $\mc
M_t = \mc M_t (a, \mb H)$ the mean one exponential martingale whose
logarithm is given by
\begin{eqnarray*}
\!\!\!\!\!\!\!\!\!\!\!\!\!\! &&
\frac a{N} \sum_{j=1}^d\sum_{x\in\bb T_N^d} 
\int_0^t  H_j (s,x/N)  dW_s^{x,x+e_j} \\
\!\!\!\!\!\!\!\!\!\!\!\!\!\! && \quad 
- N^2 \sum_{j=1}^d\sum_{x\in\bb T_N^d} \int_0^t ds\, 
\eta_s(x) [1-\eta_s(x+e_j)]  
\Big\{ e^{ a N^{-1} H_j (s,x/N)}-1 \Big\}   \\
\!\!\!\!\!\!\!\!\!\!\!\!\!\! && \qquad 
- N^2 \sum_{j=1}^d\sum_{x\in\bb T_N^d} \int_0^t ds\, 
\eta_s(x+e_j) [1-\eta_s(x)] 
\Big\{ e^{ - a N^{-1} H_j (s,x/N)} - 1 \Big\} \;.
\end{eqnarray*}
Since $H_j$ are $C^1$ functions, a Taylor expansion and a summation by
parts show that the expressions inside the integrals in the last two
terms are bounded by $C_{\mb H} a(1+a) N^d$, where $C_{\mb H}$ is a
finite constant depending only on $\mb H$. Therefore, by multiplying
by $a N^d$, adding and subtracting the appropriate integrals and
exponentiating, we get that
\begin{eqnarray*}
\!\!\!\!\!\!\!\!\!\!\!\!\!\! &&
\bb P^N_{\mu^N} \Big[ \sup_{t_k\le t\le t_{k+1}}  <\mb W_t^N, \mb H> -
<\mb W_{t_k}^N, \mb H> \; >\; \epsilon/3 \Big ] \\
\!\!\!\!\!\!\!\!\!\!\!\!\!\! && \quad
\le\; \bb P^N_{\mu^N} \Big[ \sup_{t_k\le t\le t_{k+1}} 
\mc M_t / \mc M_{t_k}  \; >\; \exp\{(1/6) a N^d \epsilon \}\Big ]
\end{eqnarray*}
provided $C_{\mb H} (1+a)\delta \le \epsilon/6$. Since $\mc M_t / \mc
M_{t_k}$ is a positive martingale equal to $1$ at time $t_{k}$, by
Doob's inequality, last expression is bounded by $\exp\{ - (1/6) a N^d
\epsilon \}$. Therefore, \eqref{f12} is less than or equal to 
$- a \epsilon/6$ for all $\delta$ small enough. This shows that the
second expression in the statement of the lemma is bounded by $- a
\epsilon/6$. Letting $a\uparrow\infty$, we conclude the proof.
\end{proof}

Standard arguments, presented in \cite[\S~10.4]{kl}, together
with Lemma \ref{s6} permit to extend the upper bound for compact sets
to closed sets.

\medskip
We conclude this section proving that we may set $J (\mb W,
\pi)=+\infty$ on the set of paths $(\mb W,\pi)$ which do not belong to
$\mf A_\gamma$.

\begin{lemma}
\label{s4}
Fix a sequence of measures $\{ \mu^N : N\ge 1\}$ and a function
$H:[0,T]\times \bb T^d \to \bb R$ in $C^{1,2}$.  Let
\begin{eqnarray}
\label{f05}
\!\!\!\!\!\!\!\!\!\!\!\!\! &&
L_T(\pi^N, H) \;=\; \< \pi^N_T, H_T\> \;-\; \< \pi^N_0, H_0\> 
\;-\; \int_0^T dt\, \< \pi^N_t, \partial_t H_t \> \;, \\
\!\!\!\!\!\!\!\!\!\!\!\!\! && 
V_T(\pi^N, \mb W^N, H) \;=\; L_T(\pi^N, H) \;-\;
\frac 1{N^{d+1}} \sum_{j=1}^d\sum_{x\in\bb T_N^d} \int_0^T
(\partial_{u_j} H) (t,x/N)  dW_t^{x,x+e_j}\;.
\nonumber
\end{eqnarray}
Then, for any $\delta>0$,
\begin{equation*}
\limsup_{N\to\infty} \frac 1{N^d}  
\log \bb P^N_{\mu^N} \Big[ \, \Big\vert V_T(\pi, \mb W^N, H) \Big\vert >
\delta \Big] \; = \; -\infty\;.
\end{equation*}
\end{lemma}

\begin{proof}
Fix a function $H:[0,T]\times \bb T^d \to \bb R$ in $C^{1,2}$. A
summation by parts shows that
\begin{eqnarray*}
\!\!\!\!\!\!\!\!\!\!\!\!\! &&
\frac 1{N^{d+1}} \sum_{j=1}^d\sum_{x\in\bb T_N^d} \int_0^T 
N\{ H(t,x+e_j/N) - H(t,x/N)\} dW_t^{x,x+e_j} \\
\!\!\!\!\!\!\!\!\!\!\!\!\! && \quad
\;=\; \frac 1{N^{d}} \sum_{x\in\bb T_N^d} \int_0^T  
H(t,x/N) \sum_{j=1}^d d\{ W_t^{x-e_j,x} - W_t^{x,x+e_j}\}\;.
\end{eqnarray*}
Since $\sum_{j=1}^d \{ W_s^{x-e_j,x} - W_s^{x,x+e_j}\}$ increases by
one each time a particle jumps to $x$ and decreases by one each time
a particle leaves $x$, this sum is equal to $\eta_s(x) -
\eta_0(x)$. In particular, an integration by parts, gives that the
previous integral is equal to $L_T(\pi^H,H)$ defined in \eqref{f05}.
Therefore, by a second order Taylor expansion,
\begin{eqnarray*}
V_T(\pi^N, \mb W^N, H) & \!\!\! = \!\!\! & 
\frac 1{N^{d+2}} \sum_{j=1}^d\sum_{x\in\bb T_N^d} 
\int_0^T (\partial^2_{u_j} H) (t,x/N)  dW_t^{x,x+e_j} \\
& \!\!\! + \!\!\! & 
\frac {o_H(1)} {N^{d+2}} \sum_{j=1}^d\sum_{x\in\bb T_N^d} 
\{ J_t^{x,x+e_j} + J_t^{x+e_j,x}\} \;,
\end{eqnarray*}
where $o_H(1)$ depends on $H$ and vanishes as $N\uparrow\infty$.

We prove that the first expression on the right hand side is
super-exponentially small, the argument for the second one being
similar. To keep notation simple, let $F_j = \partial^2_{u_j} H$. By
Chebyshev inequality, for every $a>0$,
\begin{eqnarray}
\label{f04}
\!\!\!\!\!\!\!\!\!\!\!\!\!\! &&
\bb P^N_{\mu^N} \Big[ \, \Big\vert  \frac 1{N^{d+2}} 
\sum_{j=1}^d\sum_{x\in\bb T_N^d} \int_0^T  F_j 
(t,x/N)  dW_t^{x,x+e_j}\Big\vert > \delta \Big] \;\le\; \\
\!\!\!\!\!\!\!\!\!\!\!\!\!\! && \quad
e^{-a \delta N^d} 
\bb E^N_{\mu^N} \Big[ \exp\Big\{a N^{-2} \Big\vert   
\sum_{j=1}^d\sum_{x\in\bb T_N^d} \int_0^T  F_j
(t,x/N)  dW_t^{x,x+e_j}\Big\vert \Big\} \Big]\;.
\nonumber
\end{eqnarray}
Since $e^{|x|} \le e^x + e^{-x}$, we estimate this last expectation
without the absolute value. Denote by $\bb M_T$ the mean one
exponential martingale whose logarithm is given by
\begin{eqnarray*}
\!\!\!\!\!\!\!\!\!\!\!\!\!\! &&
\frac a{N^{2}} \sum_{j=1}^d\sum_{x\in\bb T_N^d} 
\int_0^T  F_j (t,x/N)  dW_t^{x,x+e_j} \\
\!\!\!\!\!\!\!\!\!\!\!\!\!\! && \quad 
- N^2 \sum_{j=1}^d\sum_{x\in\bb T_N^d} \int_0^T dt\, 
\eta_t(x) [1-\eta_t(x+e_j)]  
\Big\{ e^{ a N^{-2} F_j (t,x/N)}-1 \Big\}   \\
\!\!\!\!\!\!\!\!\!\!\!\!\!\! && \qquad 
- N^2 \sum_{j=1}^d\sum_{x\in\bb T_N^d} \int_0^T dt\, 
\eta_t(x+e_j) [1-\eta_t(x)] 
\Big\{ e^{ - a N^{-2} F_j (t,x/N)} - 1 \Big\} \;.
\end{eqnarray*}
Since $F_j$ are continuous functions, a Taylor expansion and a
summation by parts show that the last two integrals can be written as
\begin{eqnarray*}
\Big\{ a o_F (1) + \frac{C(\mb F) a^2}{N^2} \Big\}
\sum_{j=1}^d\sum_{x\in\bb T_N^d} \int_0^T dt\, \eta_t(x)
\;\le\; \Big\{ a o_{\mb F} (1) + \frac{C(\mb F) a^2}{N^2} \Big\} d N^d T\;,
\end{eqnarray*}
where $o_{\mb F} (1)$ is an expression depending on $\mb F$ which
vanishes as $N\uparrow\infty$.

Since $\bb M_T$ is a mean one exponential martingale, the right hand
side of \eqref{f04} is bounded above by
$$
\exp a N^d \Big\{ - \delta  + o_{\mb F} (1) d  T + C(\mb F) a d N^{-2} T
\Big\}\;.
$$
In particular,
\begin{equation*}
\limsup_{N\to\infty} \frac 1{N^d}  
\log \bb P^N_{\mu^N} \Big[ \, \Big\vert  \frac 1{N^{d+2}} 
\sum_{j=1}^d\sum_{x\in\bb T_N^d} \int_0^T dt\, F_j 
(t,x/N)  dW_t^{x,x+e_j}\Big\vert > \delta \Big] \;\le\;- a \delta
\end{equation*}
for every $a>0$. This proves the lemma.
\end{proof}

From this lemma, \eqref{f13}, and Lemma~\ref{s6} it follows 
that for every closed set $F$, every $\delta>0$ and every finite
family $\{H_j ,\, 1\le j\le \ell\}$ of functions in $C^{1,2}$
\begin{equation*}
\limsup_{N\to\infty} \frac 1{N^d}  
\log \bb P^N_{\mu^N} \Big[\, (\mb W^N, \pi^N) \in F \Big] 
\;\le\; -\inf_{(\mb W,  \pi) \in F\cap \mf A_\ell} J 
(\mb W, \pi) \; ,
\end{equation*}
where 
\begin{equation*}
\mf A_\ell \;=\; \bigcap_{j=1}^\ell \Big\{ (\mb W,  \pi) \, : 
\big| V_T(\pi, \mb W, H_j) \big| \le \delta\Big\}\;.
\end{equation*}
Since this inequality holds for every $\delta>0$ and every finite
sequence $H_j$, letting $\delta\downarrow 0$ and considering a dense
family of functions $H_j$, we obtain that 
\begin{equation*}
\limsup_{N\to\infty} \frac 1{N^d}  
\log \bb P^N_{\mu^N} \Big[\, (\mb W^N, \pi^N) \in F \Big] 
\;\le\; -\inf_{(\mb W,  \pi) \in F\cap \mf A} J (\mb W, \pi) \; ,
\end{equation*}
where $\mf A$ is the set of paths $(\mb W, \pi)$ such that $\dot \pi_t
+ \nabla\cdot \dot {\mb W}_t =0$.  Up to this point, we did not need
any assumption on the sequence of initial measures $\mu^N$; but the
hypothesis that we are starting from a deterministic profile plays now
a role to replace the set $\mf A$ in the previous formula by the set
$\mf A_\gamma$, proving the upper bound of the large deviations
principle.

\subsection{The rate function.}
\label{sec4.1}
To prove the lower bound of the large deviations principle in
Theorem~\ref{s5}, we first obtain an explicit representation 
of the functional $I$ on the paths with finite rate function.

Given a path $\pi\in D([0,T];\mc F)$, we denote by $L^2(\pi)$
the Hilbert space of vector-valued functions $\mb G :[0,T]\times \bb
T^d \to \bb R^d$ endowed with the inner product $\< \cdot, \cdot
\>_\pi$ defined by
\begin{equation*}
\< \mb H, \mb G \>_\pi \;=\; \int_0^T dt \int_{\bb T^d} du\,
\chi(\pi(t,u)) \, \mb H(t,u) \cdot \mb G(t,u)\;.
\end{equation*}

Fix a pair $(\mb W,\pi)$ such that $I (\mb W, \pi) <\infty$.  In
particular $(\mb W,\pi) \in C([0,T]; \mc M_d \times \mc F)$.
Following the arguments in \cite[\S 10.5]{kl}, from Riesz
representation theorem, we derive the existence
of a function $\mb G$ in $L^2(\pi)$ such that
\begin{equation}
\label{f06}
\left\{
\begin{array}{l}
\vphantom{\Big\{}
I (\mb W, \pi) \;=\; \frac 12 \int_0^T dt 
\, \< \chi(\pi_t), | \mb G_t| ^2 \>   \\
\vphantom{\Big\{}
\partial_t \mb W_t \;+\; (1/2) \nabla \pi_t \;=\; \chi(\pi_t) \mb G_t\;.
\end{array}
\right.
\end{equation}
where the last equation has to be understood in the weak sense: 
for each $\mb H \in C^1(\bb T^d;\bb R^d)$ and each $0\le s\le t\le T$,
we have 
\begin{equation*}
\<\mb W_t, \mb H\> \;-\; \<\mb W_s, \mb H\>  \;=\;
(1/2) \int_s^t dr\, \<\pi_r , \nabla \cdot \mb H\> 
\; +\;  \int_s^t dr\, \< \chi(\pi_r) , \mb G_r \cdot \mb H\>\;.
\end{equation*}

\subsection{The lower bound}
\label{sec4}
In this Subsection we prove the lower bound in Theorem~\ref{s5}.
Denote by $\mc S$ the set of trajectories $(\mb W, \pi)$ in $\mf
A_\gamma$ for which there exists a vector-valued function $\mb G$ in
$C^{1,1}([0,T]\times \bb T^d)$ such that $(\mb W, \pi)$ is the
solution of \eqref{f06}. 

For paths $(\mb W, \pi)$ in $\mc S$, we may repeat the arguments
presented in the proof of the lower bound in the large deviations
principle for the empirical density in \cite[\S~10.5]{kl} to
conclude that for each open set $G$
\begin{equation*}
\liminf_{N\to\infty} \frac 1{N^d}  
\log \bb P^N_{\eta^N} \Big[\, (\mb W^N, \pi^N) \in G \Big] 
\;\ge\; - \inf_{(\mb W, \pi) \in G \cap \mc S} I (\mb W,
\pi)\;. 
\end{equation*}

To conclude the proof, it remains to show that for all pairs $(\mb W,
\pi)$ with finite rate function, $I (\mb W, \pi) < \infty$, there
exists a sequence $(\mb W_k, \pi_k)$ in $\mc S$ converging to $(\mb W,
\pi)$ and such that $\lim_{k\to\infty} I (\mb W_k, \pi_k) = I (\mb W,
\pi)$. When this occurs, we shall say that the sequence $(\mb W_k,
\pi_k)$ $I$-converges to $(\mb W, \pi)$.
In the context of the symmetric exclusion process, the argument is not
too difficult because the rate function is convex. We follow the proof
of the lower bound presented in \cite{bdgjl3}.

The proof is divided in two steps. We first show that there exists
a sequence $(\mb W_k, \pi_k)$ which $I$-converges to $(\mb W, \pi)$
and such that, for each $k$, $\pi_k$ is bounded away from $0$ and $1$
uniformly in $[0,T]\times \bb T^d$. To do it, following \cite{bdgjl3},
we consider a convex combination of $(\mb W, \pi)$ with the solution
of the hydrodynamic equation \eqref{f03} with external field $\mb F=0$
and initial condition $\pi_0=\gamma$, $\mb W_0 = 0$.

Consider now a pair $(\mb W, \pi)$ whose empirical density $\pi$ is
bounded away from $0$ and $1$. Since $I (\mb W, \pi)$ is finite, by
Subsection~\ref{sec4.1}, there exists a vector-valued function
$\mb G$ in $L^2(\pi)$ satisfying \eqref{f06}.  Since $\pi$ is bounded
away from $0$ and $1$, $L^2(\pi)$ coincides with the usual $L^2$ space
associated to the Lebesgue measure on $[0,T]\times \bb T^d$.  Consider
a sequence of smooth vector-valued functions $\mb G_n : [0,T]\times
\bb T^d \to \bb R^d$ converging in $L^2$ to $\mb G$ and denote by
$(\mb W^n, \pi^n)$ the pair in $\mf A_\gamma$ which solves \eqref{f06}
with $\mb G_n$ instead of $\mb G$. Repeating the arguments presented
in \cite[\S~3.6]{bdgjl3}, one can prove that $(\mb W^n,
\pi^n)$ $I$-converges to $(\mb W, \pi)$. This concludes the proof of
the lower bound.

\subsection{Large deviations for the empirical density.}

In this subsection we show that the large deviation principle for the
empirical density, proven in \cite{kov}, follows from Theorem~\ref{s5}.
Indeed, the large deviations principle for the empirical density 
can be recovered from the one for the current by 
the contraction principle. The rate function 
$\mc I$ is given by the variational formula
\begin{equation}
\label{f08}
\mc I (\pi) \;=\; \inf_{\mb W \in \mf W_\pi} 
I (\mb W, \pi)\;,
\end{equation}
where $\mf W_\pi$ stands for set of currents $\mb W$ satisfying
$\dot \pi_t + \nabla\cdot \dot {\mb W}_t =0$, as formulated in
\eqref{f09}. 

This variational problem is simple to solve. Let us first assume that
$\pi$ is smooth and bounded away from $0$ and $1$. Fix a current $\mb
W$ in $\mf W_\pi$ and denote by $\mb G$ the external field associated
to $\mb W$ through equation \eqref{f06}. For $0\le t\le T$, let $H_t$
be the solution of the elliptic equation
\begin{equation*}
\nabla \cdot \big( \chi(\pi_t) \mb G_t \big) \;=\; 
\nabla \cdot \big( \chi(\pi_t) \nabla H_t \big)
\end{equation*}
and set $\mb F_t = \chi(\pi_t) \{ \mb G_t - \nabla H_t\}$. By
definition $\nabla \cdot \mb F_t=0$. Let $\mb w$ be the current
defined by
\begin{equation}
\label{f07}
\dot {\mb w}_t \;+\; (1/2) \nabla \pi_t \;=\; \chi(\pi_t) \nabla H_t\;.
\end{equation}
$\mb w$ belongs to $\mf W_\pi$ because by construction $\nabla \cdot
\dot {\mb w} = \nabla \cdot \dot {\mb W}$. Moreover, by the explicit
formula \eqref{f06} for the rate function and by definition of $\mb F$,
\begin{eqnarray*}
I (\mb W, \pi) &\!\!\! = \!\!\! &  \frac 12 \int_0^T dt 
\, \< \chi(\pi_t), | \mb G_t| ^2 \> \\
&\!\!\! = \!\!\! &  \frac 12 \int_0^T dt 
\Big\{ \< \chi(\pi_t), | \nabla H_t| ^2 \>
+ 2  \< \mb F_t \cdot  \nabla H_t \>  + 
\< \chi(\pi_t)^{-1} |\mb F_t| ^2 \> \Big\}\;.
\end{eqnarray*}
Since $\nabla \cdot \mb F=0$, an integration by parts shows that the
cross term vanishes. On the other hand, by the explicit formula
\eqref{f06} for the rate function and by \eqref{f07}, the first term
on the right hand side is $I (\mb w, \pi)$. Thus,
\begin{equation*}
I (\mb W, \pi)\; \ge \; I (\mb w, \pi)\; .
\end{equation*}
In particular, in the variational problem \eqref{f08}, we can restrict
our attention to currents $\mb W$ for which the associated external
fields $\mb G$ are in gradient form.

Now, consider two currents $\mb W^1$, $\mb W^2$ in $\mf W_\pi$
and assume that both external fields $\mb G_1$, $\mb G_2$ associated
to these currents through \eqref{f06} are in gradient form: $\mb G_j
= \nabla H^j$, $j=1$, $2$. Taking the divergence of \eqref{f06} and
recalling that $\dot \pi_t + \nabla\cdot \dot {\mb W}_t =0$, we obtain
that
\begin{equation*}
\dot \pi_t \;=\; (1/2) \Delta \pi_t - 
\nabla \cdot \big( \chi(\pi_t) \nabla H^j_t \big) 
\end{equation*}
for $j=1$, $2$ and each $0\le t\le T$. In particular, 
$\nabla \cdot \big( \chi(\pi_t) \nabla [H^1_t - H^2_t]  \big)=0$. 
Taking the inner product with
respect to $H^1_t-H^2_t$ and integrating by parts, we get that
\begin{equation*}
\int_{\bb T^d} du\, \chi(\pi_t) \vert \nabla H^1_t -
\nabla H^2_t\vert^2 \;=\; 0
\end{equation*}
for every $0\le t\le T$.  In particular, $I (\mb W^1, \pi) =
I (\mb W^2, \pi)$. This proves that the variational problem
\eqref{f08} is attained on currents for which the associated external
field is in gradient form:
\begin{equation*}
\mc I (\pi) \;=\; \inf_{\mb W \in \mf W_\pi} \: I (\mb W, \pi)
\;=\; \frac 12 \int_0^T dt \< \chi(\pi_t) \vert \nabla H_t\vert^2 \>
\;,
\end{equation*}
where $H_t$ is given by
\begin{equation*}
\dot \pi_t \;=\; (1/2) \Delta \pi_t 
- \nabla \cdot \big( \chi(\pi_t)  \nabla H_t \big) \;.
\end{equation*}
This is exactly the large deviations rate function for the empirical
density obtained in \cite{kov}. This identity has been obtained for
smooth paths $\pi$ bounded away from 0 and 1. 
However, by the arguments of the previous subsection, 
we can extend it to all paths $\pi$.

\section{Large deviations of the mean current on a 
long time interval}    
\label{sec:6}

In this Section we investigate the large deviations properties of the 
mean empirical current $\mb W^N_T/T$ as we let \emph{first} $N\to\infty$
and \emph{then} $T\to\infty$. We emphasize that the analysis carried
out in this section does not depend on the details of the symmetric
simple exclusion process so that it holds in a general setting.

Given a profile $\gamma \in C^2(\bb T^d)$, $T>0$, and 
$\mb W\in D\big([0,T];\mc M_d\big)$, let $\pi\in D\big([0,T];\mc F\big)$ be the
solution of \eqref{f09} and  denote by
 $I_{[0,T]}(\mb W|\gamma)$ the functional defined in \eqref{f10}, 
in which we made explicit the dependence on the time interval $[0,T]$ 
and on the initial profile $\gamma$.  
We define $\Phi_T (\cdot |\gamma) : \mc M_d\to [0,+\infty]$ as
the functional
\begin{equation}
\label{IT}
\Phi_T(\mb J|\gamma) = 
\frac 1T \inf_{\mb W \in \mc A_{T,\mb J} } I_{[0,T]}  (\mb W | \gamma)
\end{equation}
where
$$
\mc A_{T,\mb J} := 
\Big\{  \mb W \in D\big([0,T];\mc M_d\big) \,:\, \mb W_T = T \mb J  \Big\}
$$
Recalling  that the set $\mc B$ of divergence free measures has been
defined in \eqref{dB} we also define
\begin{equation}
\label{tilphi}
\tilde \Phi (\mb J |\gamma) :=
\begin{cases}
{\displaystyle \inf_{T>0} \Phi_T(\mb J|\gamma)} & \text{if $\mb J \in \mc B$}
\\
 +\infty & \text{otherwise} 
\end{cases}
\end{equation}
Denote finally by  $\Phi (\mb J |\gamma) := \sup_{U \ni \mb J} \,\,
\inf_{\mb{J}'\in U}\,\, \tilde \Phi (\mb{J}' |\gamma)$, 
where $U\subset \mc M_d$ is open, the lower semi-continuous envelope 
of $\tilde \Phi (\cdot |\gamma)$.

\begin{remark}
\label{t:gindip}
The functional $\Phi(\cdot|\gamma)$ depends 
on the initial condition $\gamma$ only through its total mass 
$m=\int\!du \, \gamma(u)$. This holds in the present setting of
periodic boundary conditions; in the case of Dirichlet boundary
conditions, when the density is fixed at the boundary,
$\Phi(\cdot|\gamma)$ would be completely independent on $\gamma$. 
Furthermore the functional $\Phi(\cdot|\gamma)$ is convex. 
\end{remark}

\begin{theorem}
\label{t:LDphi}
Let $\gamma\in C^2 (\bb T^d)$ bounded away from $0$ and $1$, and 
$\eta^N\in \mc X_N$  a sequence such that $\pi^N(\eta^N)\to
\gamma$ in $\mc F$.
Then, for each closed set $C$ and every open set $U$ of $\mc M_d$,
we have
\begin{eqnarray*}
&&
\varlimsup_{T\to\infty} \varlimsup_{N\to\infty} 
\frac 1{T N^d} \log 
\bb P_{\eta^N}^N
\Big[ \frac 1T {\bf W}^N_T \in C \Big] \;\le\; 
- \inf_{{\mb J}\in C} \Phi ({\mb J}|\gamma) \;, 
\\ 
&&\qquad
\varliminf_{T\to\infty} \varliminf_{N\to\infty} 
\frac 1{T N^d} \log \bb P_{\eta^N}^N
\Big[ \frac 1T {\mb W}^N_T \in U \Big] \;\ge\; 
- \inf_{{\mb J}\in U} \Phi({\mb J}|\gamma) \;. 
\end{eqnarray*}
\end{theorem}

If $T$ were fixed, the above large deviation principle would directly 
follow from Theorem~\ref{s5}. The asymptotic $T\to\infty$ is related 
to the so-called $\Gamma$-convergence of the rate functions. 
We first discuss this issue in a general setting. 
Let $\mc X$ be a metric space, we recall that a sequence $F_T: \mc X
\to  [0,+\infty]$ of functions  $\Gamma$-converges to 
$F: \mc X \to [0,+\infty]$ as $T\to\infty$ 
iff for each $x\in\mc X$ the following holds
\begin{eqnarray}
\label{Gub}
\!\!\!\!\!\!\!
&\text{for any sequence $x_T\to x$ we have} 
& \;\;
F(x) \le \varliminf_{T\to\infty} F_T(x_T) \quad\quad \quad\quad 
\\
\label{Glb}
\!\!\!\!\!\!\!
& \text{there exists a sequence $x_T\to x$ such that} 
& \;\;
F(x) \ge \varlimsup_{T\to\infty} F_T(x_T)  
\quad\quad \quad\quad 
\end{eqnarray}

\begin{lemma}
\label{t:GtLD}
Let $P_{N,T}$ be a two parameter family of probabilities on $\mc X$
endowed with its Borel $\sigma$-algebra. Assume that for each fixed 
$T>0$ the family $\{P_{N,T}\}_{N\in\bb N}$ 
satisfies the weak large deviation
principle with rate function $ T\, F_T$, that is 
for each $K$ compact and $U$ open in $\mc X$ we have
\begin{eqnarray}
\label{LDTub}
&&
\varlimsup_{N\to\infty} 
\frac 1N \log P_{N,T}(K) \le - T \, \inf_{x\in K} F_T(x)
\\
\label{LDTlb}
&& \qquad
\varliminf_{N\to\infty} 
\frac 1N \log P_{N,T}(U) \ge - T \, \inf_{x\in U} F_T(x)
\quad\quad \quad\quad 
\end{eqnarray}
Assume also that the sequence $F_T$ $\,\Gamma$-converges to $F$ 
as $T\to\infty$. Then for each $K$ compact and $U$ open in
$\mc X$ we have
\begin{eqnarray}
\label{LDub}
&& 
\varlimsup_{T\to\infty} \varlimsup_{N\to\infty} 
\frac 1{NT} \log P_{N,T}(K) \le - \inf_{x\in K} F(x)
\\
\label{LDlb}
&& \qquad
\varliminf_{T\to\infty} \varliminf_{N\to\infty} 
\frac 1{NT} \log P_{N,T}(U) \ge - \inf_{x\in U} F(x)
\end{eqnarray}
\end{lemma}

\begin{proof}
To deduce \eqref{LDub}--\eqref{LDlb} from \eqref{LDTub}--\eqref{LDTlb}
we need to show that for each $K$ compact and $U$ open in $\mc X$ we have
\begin{eqnarray*}
&& 
\varliminf_{T\to\infty}  \inf_{x\in K} F_T(x) 
\ge \inf_{x\in K} F(x)
\\
&& \qquad
\varlimsup_{T\to\infty} 
 \inf_{x\in U} F_T(x) \le \inf_{x\in U} F(x)
\end{eqnarray*}
These bounds are a direct consequence of \eqref{Gub}--\eqref{Glb},
see e.g.\ \cite[Prop.~1.18]{B}. 
\end{proof}

The following Lemma follows from Theorem~\ref{s5} by contraction
principle.

\begin{lemma}
\label{t:LDT}
Consider a profile $\gamma\in C^2(\bb T^d)$ bounded away from $0$ and $1$ 
and a sequence $\{\eta^N : N\ge 1\}$ such that $\pi^N(\eta^N)\to
\gamma$ in $\mc F$. Then for each $T>0$ we have the following 
large deviation principle for the mean empirical current $\mb W^N_T/T$.
For each $C$ closed and each $U$ open in $\mc M_d$
\begin{eqnarray*}
&&
\varlimsup_{N\to\infty} \frac 1{N^d} 
\log \bb P_{\eta^N}^N \Big[ \frac{\mb W^N_T}{T} \in C \Big] 
\le - T \, \inf_{\mb J\in C} \Phi_T(\mb J|\gamma)
\\
&&\qquad
\varliminf_{N\to\infty} \frac 1{N^d} 
\log \bb P_{\eta^N}^N \Big[ \frac{\mb W^N_T}{T} \in U \Big] 
\ge - T \, \inf_{\mb J\in U} \Phi_T(\mb J|\gamma)
\end{eqnarray*}
where we recall the functional $\Phi_T$ is defined in \eqref{IT}.
\end{lemma}

\begin{proposition}
\label{t:Gcon}
Let $\gamma\in C^2 (\bb T^d)$ bounded away from $0$ and $1$.
The sequence of functionals
$\Phi_T(\cdot|\gamma)$ $\,\,\Gamma$-converges to the functional 
$\Phi(\cdot|\gamma)$ defined after \eqref{tilphi}.
\end{proposition}

The previous Proposition, together with Lemmata~\ref{t:GtLD} 
and \ref{t:LDT}, 
proves the ``large deviations principle'' stated in
Theorem~\ref{t:LDphi} for compact sets. 
For its proof, we need a few preliminary Lemmata.
For each $\mb J \not\in \mc B$, since the empirical density is
bounded, by the continuity equation \eqref{f09}, 
we have that $\Phi_T(\mb J|\gamma) = +\infty$ if $T$ is sufficiently
large. We next show that this holds uniformly for all $\mb J$
whose distance from $\mc B$ is uniformly bounded below. 
To this end we introduce the following metric on $\mc M_d$. 
Pick a sequence of smooth vector fields 
${\mb G}_k\in C^1\big(\bb T^d;\bb R^d\big)$, $k\ge 1$, 
dense in the unit ball of $C\big(\bb T^d;\bb R^d\big)$; 
for $\mb J,\mb{J}' \in \mc M_d$ we then define
$$
\varrho(\mb J,\mb J') = \sum_{k=1}^\infty 
\frac 1{2^k}   \: 1\wedge \big| \langle \mb J-\mb J', \mb G_k \rangle \big|
$$
It is easy to show that $\varrho$ is a metric inducing 
the weak topology of $\mc M_d$.

\begin{lemma}
\label{t:d>d}
For each $\delta\in (0,1)$ 
there exists $T_0=T_0(\delta)\in \bb R_+$ such that for any
$T\ge T_0$ we have $\Phi_T(\mb J|\gamma) = +\infty$ 
for any $\mb J\in \mc M_d$ such that $\varrho(\mb J, \mc B) \ge \delta$. 
\end{lemma}

\begin{proof}

Given $\mb J\in \mc M_d$ let us denote by $\widehat{\mb J}\in \mc B$ its
projection on the subspace $\mc B$. If 
$\mb J(du) = \mb j \,du$ for some $\mb j \in L_2(\bb T^d;\bb R^d)$ 
this is simply the orthogonal projection of $\mb j$ to $\mc B$. 
In general $\widehat{\mb J}\in \mc M_d$  is defined by 
$\langle \widehat{\mb J}, \nabla F \rangle = 0$ for any
$F\in C^1\big(\bb T^d\big)$ and $\langle \widehat{\mb J}, \mb E \rangle 
= \langle \mb J, \mb E \rangle$ for any 
$\mb E\in C^1(\bb T^d;\bb R^d)$ such that 
$\nabla\cdot \mb E =0$. It is easy to verify that
$\widehat{\mb J}$ is uniquely defined by the above requirements. 
We then have
$$
\delta \le \varrho(\mb J, \mc B) \le  
\varrho(\mb J,  \widehat{\mb J}) 
= \sum_{k=1}^\infty \frac 1 {2^k}\: 
1\wedge \big| \langle \mb J - \widehat{\mb J}, \mb G_k \rangle \big| 
= \sum_{k=1}^\infty \frac 1 {2^k} \:
1\wedge \big| \langle \mb J, \nabla F_k \rangle \big|  
$$
where $F_k \in C^1\big(\bb T^d\big)$ is obtained from $\mb G_k$ by
solving the Poisson equation $\Delta F_k = \nabla \cdot \mb G_k$ so
that $\mb G_k = \nabla F_k + \mb E_k$ with $\nabla \cdot {\mb E}_k = 0$.  
Note that $\|F_k\|_{L_2} \le C_0 \| \mb G_k \|_{L_2}\le C_0$ for some
constant $C_0$ not depending on $k$. 

From the previous inequality we get that there exists  
$\bar k = \bar k (\mb J)$ such that  
$\big| \langle \mb J, \nabla F_{\bar k} \rangle \big| \ge \delta$. 
Let $\mb W\in \mc A_{T,\mb J}$ and denote by $\pi_t$, $t\in [0,T]$ the
corresponding solution of  the continuity equation \eqref{f09}.
By choosing $F=F_{\bar k}$ and using that $\mb W_T = T \mb J$, from
\eqref{f09} we get 
$$
\langle \pi_T- \gamma, F_{\bar k} \rangle = 
T\, \langle \mb J , \nabla F_{\bar k} \rangle  
$$
Since $-1\le\pi_T- \gamma \le 1$, the absolute value of the l.h.s.\ above
is bounded above by $C_0$. On the other hand, since 
$\varrho(\mb J,\mc B)\ge \delta$, the absolute value of the r.h.s.\ above 
is bounded below by $\delta T$. 
By taking $T_0 > C_0 \, \delta^{-1}$, the Lemma follows.
\end{proof}

\begin{lemma}
\label{t:l2b}
Consider a profile $\gamma\in C^2(\bb T^d)$ bounded away from $0$ and $1$  
and let $\mb J \in \mc B$, $T>0$.
We then have $\Phi_T( \mb J|\gamma) < +\infty $ iff $\mb J(du)=\mb j\, du$ 
for some  $\mb j \in L_2\big(\bb T^d; \bb R^d\big)$. Moreover there exists a
constant $C_1 \in (0,\infty)$ (depending on $\gamma$) such that for any $T>0$
and any $\mb J(du)= \mb j\, du$ we have 
\begin{equation}
\label{el2}
\frac 1{C_1} \langle \mb j, \mb j\rangle  
\le \Phi_T ( \mb J | \gamma) 
\le C_1   \big[ \langle \mb j, \mb j\rangle +  1 \big]
\end{equation}
\end{lemma}

\begin{proof} 
Let $\mb W\in \mc A_{T,\mb J}$; by choosing the vector field $\mb F$ in 
the variation expression \eqref{f10} constant in time and divergence
free we get 
$$
\frac 1T  I_{[0,T]} (\mb W|\gamma) 
\ge  \langle \mb J, \mb F \rangle - 
\frac 12    \, \frac 14 \langle \mb F, \mb F \rangle
$$
where we used that $\chi(\pi) \le 1/4$. Recalling that $\mb J\in \mc B$, 
by optimizing over $\mb F\in \mc B$ we see that $\Phi_T(\mb J|\gamma)
=+\infty$ unless $\mb J (du)= \mb j \, du$ for some 
$\mb j \in L_2\big(\bb T^d; \bb R^d\big)$. In
fact this argument also proves the first inequality in
\eqref{el2}. To prove the second inequality in \eqref{el2} it is
enough to construct an appropriate path $\mb W\in \mc A_{T,\mb J}$; we
simply take $\mb W_t(du)= t\, \mb j\, du$. 
The solution of the continuity equation
\eqref{f09} is then given by $\pi_t=\gamma$ and, by \eqref{f06}, 
$$
\frac 1T I_{[0,T]}(\mb W |\gamma) = 
\frac 12 
\Big\langle  \frac 1 {\chi(\gamma)}, 
\Big| \mb j + \frac 12 \nabla \gamma \Big|^2 \Big\rangle
$$
Recalling that $\gamma$ is bounded away from $0$ and $1$,
$\chi(\gamma)=\gamma(1-\gamma)$,  the result follows. 
\end{proof}

We next show that on divergence free measures $\mb J$ 
the functional $T \, \Phi_T(\cdot|\gamma)$ is sub-additive.

\begin{lemma}
\label{t:subadd}
Consider a profile $\gamma\in C^2(\bb T^d)$ bounded away from $0$ and
$1$ and let $\mb J \in \mc B$. Then, for each $T,S>0$, we have
$$
(T+S) \, \Phi_{T+S} (\mb J|\gamma) 
\le T \, \Phi_{T} (\mb J|\gamma) + S \, \Phi_{S} (\mb J|\gamma)
$$
\end{lemma}

\begin{proof}
By \eqref{IT},  
given $\varepsilon>0$ there are $\mb W^1 \in \mc A_{T,\mb J}$ and 
$\mb W^2 \in \mc A_{S,\mb J}$ such that 
$$
\Phi_{T} (\mb J |\gamma) \ge  T \, I_{[0,T]}(\mb W^1|\gamma) 
- \frac 12 \, \varepsilon
\qquad
\Phi_{S} (\mb J |\gamma) \ge S\, I_{[0,S]}(\mb W^2|\gamma) -\frac 12\, 
\varepsilon
$$
Let $\mb W_t$, $t\in [0,T+S]$ the path obtained by gluing $\mb W^1$ with
$\mb W^2$, i.e.\ we define
$\mb W_t := \mb W^1_{t\wedge T} + \id_{[T,T+S]}(t) \mb W^2_{t-T}$,
and denote by $\pi_t$ the corresponding solution of the
continuity equation \eqref{f09}. Then $\mb W\in \mc A_{T+S,\mb J}$ and, by the 
invariance of $I_{[0,T]}$ with respect to time shifts, 
\begin{eqnarray*}
(T+S)\, \Phi_{T+S} (\mb J|\gamma) 
&\le &  I_{[0,T+S]}(\mb W|\gamma) 
\;= \;  I_{[0,T]}(\mb W^1|\gamma)  + I_{[0,S]}(\mb W^2|\pi_T) 
\\
&\le &
T \, \Phi_{T} (\mb J|\gamma) + S\, \Phi_{S} (\mb J|\gamma) 
+ \varepsilon
\end{eqnarray*} 
where we used that $\mb J\in \mc B$, $\mb W^1\in \mc A_{T,\mb J}$
implies $\pi_T= \gamma$. 
\end{proof}

Recall that $\tilde \Phi (\cdot |\gamma)$ is defined in \eqref{tilphi} 
and note that, by Lemmata~ \ref{t:d>d} and \ref{t:l2b}, 
$\tilde \Phi (\mb J|\gamma)$ equals 
$+\infty$ unless $\mb J(du) =\mb j\, du$ 
for some divergence free $\mb j\in L_2(\bb T^d;\bb R^d)$. 
By the sub-additivity proven above we have that 
$\tilde \Phi (\mb J|\gamma) = \lim_{T\to\infty} \Phi_T(\mb J|\gamma)$ 
pointwise in $\mb J$.
However the pointwise convergences does not imply the
$\Gamma$-convergence and some more efforts are required.

\begin{lemma}
\label{t:Gub}
Consider a profile $\gamma\in C^2(\bb T^d)$ bounded away from $0$ and
$1$, and let $\bar{\mb J}\in \mc B$. 
Then, for each open neighborhood $U$ of $\bar{\mb J}$ we have 
$$
\varlimsup_{T\to\infty} \inf_{\mb J\in U} \Phi_T(\mb J|\gamma) 
\le \inf_{\mb J\in U} \tilde \Phi (\mb J|\gamma)
$$
\end{lemma}

\begin{proof}
Thanks to Lemma~\ref{t:l2b} we can assume that the set $U$ is a 
bounded subset of 
$L_2(\bb T^d;\bb R^d)$.
Pick $S>0$ and let $k:=[T/S]$, $R:= T-kS \in [0,S)$; 
by Lemma~\ref{t:subadd} we get 
\begin{eqnarray*}
&&\varlimsup_{T\to\infty} \inf_{\mb J\in U} \Phi_T(\mb J|\gamma) 
\le 
\varlimsup_{k\to\infty} \inf_{\mb J\in U \cap \mc B} 
\Big[ \frac{k S}{kS+R}   \Phi_S (\mb J|\gamma) 
+  \frac{R}{kS+R}\Phi_R (\mb J|\gamma)  \Big]
\\
&&\quad\quad\quad
\le \varlimsup_{k\to\infty} \Big[ \frac{k S}{kS+R}   
\inf_{\mb J\in U \cap \mc B}  \Phi_S (\mb J|\gamma) 
+ \frac{R}{kS+R} \sup_{\mb J\in U \cap \mc B}  \Phi_R (\mb J|\gamma) 
\Big]
\\
&&\quad\quad\quad
= \inf_{\mb J\in U \cap \mc B}  \Phi_S (\mb J|\gamma) 
\end{eqnarray*}
where we used again Lemma~\ref{t:l2b} to get  
$\sup_{\mb J\in U \cap \mc B}  \Phi_R (\mb J|\gamma) < \infty$. 
By taking the infimum over $S>0$ we get the result. 
\end{proof}

\begin{lemma}
\label{t:cammino}
Given $m\in (0,1)$, let $S_m :\mc F_m \to \bb R_+$ be the
functional
$$
S_m (\rho) := \int_{\bb T^d} \!du\: 
\Big\{ \rho(u) \log \frac{\rho(u)}{m} 
+ [1-\rho(u)] \log \frac{1-\rho(u)}{1-m}  \Big\}
$$
Then for each $\delta>0$ there exists $T_0=T_0(\delta)>0$
such that the following holds. For each $\gamma_1,\gamma_2 \in \mc
F_m$ there exists a path $(\mb W,\pi)\in  
C([0,T_0];\mc M_d\times \mc F) \cap \mf A_{\gamma_1}$,
such that: $\pi_{0}=\gamma_1$, $\pi_{T_0}=\gamma_2$, 
\begin{equation}
\label{cucam}
\big|\langle {\mb W}_{T_0}, {\mb F}\rangle \big| \le  
\frac 12 T_0 \|\nabla \cdot {\mb F}\|_{L_2}
+ \delta \|{\mb F}\|_{L_2}
\qquad \text{ for any }\;  {\mb F}\in C^1(\bb T^d;\bb R^d),
\end{equation}
and 
\begin{equation}
\label{eccam}
I_{[0,T_0]}\big( ({\mb W} ,\pi) \big| \gamma_1 \big) 
\le S_m(\gamma_2) + \delta
\end{equation}
\end{lemma}

\begin{proof}
The strategy to construct the path $\pi$ is the following. Starting
from $\gamma_1$ we follow the hydrodynamic equation \eqref{f02} until we
reach a small neighborhood (in a strong topology) of the constant
profile $m$, paying no cost, then we move ``straight'', paying only a
small cost,   to a suitable point in that small neighborhood which is
chosen so that starting from it we can follow the time reversed
hydrodynamic equation to get to $\gamma_2$; the cost of this portion of
the path is $S_m(\gamma_2)$. The current ${\mb W}$ is chosen so
\eqref{f09}, \eqref{f07} hold, i.e.\ is the one whose cost is
minimal among the ones compatible with the density path $\pi$.

Let $\lambda\in \mc F_m$ and denote by $P_t \lambda$ the solution of
the Cauchy problem \eqref{f02} ($P_t$ is indeed the heat semigroup on
$\bb T^d$). By the regularizing properties of the heat semigroup,
given $\delta_1>0$ there exists a time $T_1$ such that 
$\|P_t \lambda\|_{H_1} + \|P_t \lambda - m \|_{\infty} \le \delta_1$, for any
$t\ge T_1$. 
Here $\|\varphi\|_{H_1}=\|\nabla \varphi\|_{L_2}$ is the standard Sobolev
norm on $\bb T^d$ and the time $T_1$ is independent on $\lambda$ because
$0\le \lambda \le 1$.  
We now choose $\delta_1 < (1/2)\, [m \wedge (1-m)]$ and let 
$T_0:= 2 T_1 +1$, $\bar\gamma_i =
P_{T_1} \gamma_i$, $i=1,2$. The density path $\pi$ is then constructed
as
$$
\pi_t :=
\begin{cases}
P_t \gamma_1 &  \text{for} \quad t\in [0,T_1] \\ 
\bar\gamma_1 [1-(t-T_1)] + \bar\gamma_2 [t-T_1]
&  \text{for} \quad t\in (T_1,T_1+1) \\ 
P_{T_0-t} \gamma_2  &  \text{for} \quad t\in [T_0-T_1,T_0]
\end{cases}
$$
while the associated current path ${\mb W}$ is such that 
$$
\dot{{\mb W}}_t=
\begin{cases}
-\frac 12 \nabla P_t \gamma_1 &  \text{for} \quad t\in [0,T_1] \\ 
-\frac 12 \nabla \pi_t + \chi(\pi_t) \nabla H_t   &  \text{for}
\quad t\in (T_1,T_1+1) \\ 
\frac 12 \nabla P_{T_0-t} \gamma_2  &  \text{for} \quad t\in [T_0-T_1,T_0]
\end{cases}
 $$
where $H\in C^{1,2}((T_1,T_1+1)\times\bb T^d)$ solves 
$$
\nabla\cdot \big[ \chi(\pi_t)\nabla H_t \big]= -\partial_t \pi_t +\frac 12 
\Delta \pi_t
$$
Note there exists a unique solution since the r.h.s.\ 
is orthogonal to the constants (note $\bar\gamma_i\in \mc F_m$, $i=1,2$).

It is straightforward to verify that $(\pi,{\mb W})\in \mf
A_{\gamma_1}$, i.e.\ the continuity equation \eqref{f09} holds. 
Thanks to the invariance of $I$ with respect to time shifts, 
the cost of this path is 
\begin{eqnarray}
\label{ccam}
I_{[0,T_0]}\big( ({\mb W}, \pi) \big| \gamma_1 \big)  
&=& 
I_{[0,T_1]}\big( ({\mb W} , \pi )  \big| \gamma_1 \big) 
+I_{[0,1]}\big( ({\mb W}_{\cdot-T_1}, \pi_{\cdot -T_1} )
\big| \bar\gamma_1 \big) 
\nonumber 
\\
&&
+\; I_{[0,T_1]}\big(  ({\mb W}_{\cdot - (T_1+1)}, \pi_{\cdot-(T_1+1)}) 
\big| \bar\gamma_2 \big) 
\end{eqnarray}
Since in the time interval $[0,T_1]$ the path follows the hydrodynamic
equation, the first term on the r.h.s.\ of \eqref{ccam} vanishes. By
considering the time reversal of the portion of the path in the
interval $[T_0-T_1,T_0]$, it is straightforward to verify, see
\cite[Lemma~5.4]{bdgjl3}, that
$$
 I_{[0,T_1]}\big(  ({\mb W}_{\cdot - (T_1+1)}, \pi_{\cdot-(T_1+1)})  
\big| \bar\gamma_2 \big) 
 = S_m(\gamma_2)-S_m(\bar\gamma_2) \le  S_m(\gamma_2) 
$$

It remains to bound the second term on the r.h.s.\ of
\eqref{ccam}. Note that, by construction, in the interval
$(T_1,T_1+1)$ we have 
$$
\inf_{u} \bar\gamma_1(u) \wedge \inf_{u} \bar\gamma_2(u)
\le \pi_t \le \sup_{u} \bar\gamma_1(u) \vee \sup_{u} \bar\gamma_2(u)
$$
which implies, by the choice of $\delta_1$, 
$$
\frac 12 \, \big[ m\wedge (1-m) \big] 
\le \pi_t \le 1 -  \frac 12 \, \big[ m\wedge (1-m) \big]
$$
i.e.\ the the path $\pi_t$ is uniformly bounded away from 0 and 1. By the same
computations as in \cite[Lemma~5.7]{bdgjl3}, it is not difficult to
show that there exists a constant $C>0$ depending only on $m$ such that 
\begin{eqnarray*}
I_{[0,1]}\big( {\mb W}_{\cdot-T_1}, \pi_{\cdot -T_1} 
\big| \bar\gamma_1 \big) 
&\le & C \int_0^1\!dt\: 
\big\|\partial_{t} \pi_{t+T_1} +\frac 12 \Delta\pi_{t+T_1} 
\big\|^2_{H_{-1}}
\\
&\le & C' \big[ 
\big\| \bar\gamma_2\big\|_{H_1}^2 +  \big\| \bar\gamma_1\big\|_{H_1}^2
\big]
\end{eqnarray*}
which, by taking $\delta_1$ small enough, concludes the proof of
\eqref{eccam}. The bound \eqref{cucam} follows easily from the
construction of the path $\mb W$ by using Cauchy--Schwartz and what proven
above.  
\end{proof}

\begin{lemma}
\label{t:Glb}
Consider a profile $\gamma \in C^2(\bb T^d)$ bounded away from $0$ and
$1$, and let $\bar{\mb J} \in \mc B$. 
Then, for each open neighborhood $U$ of $\bar{\mb J}$ we have 
$$
\varliminf_{T\to\infty} \inf_{\mb J\in U} \Phi_T(\mb J|\gamma) 
\ge \inf_{\mb J\in U} \tilde \Phi (\mb J|\gamma)
$$
\end{lemma}

\begin{proof}
Recalling  definition \eqref{tilphi}, it is enough to show
\begin{equation}
\label{toGlb}
\varliminf_{T\to\infty} \inf_{\mb J\in U} \Phi_T(\mb J|\gamma) 
\ge \varliminf_{T\to\infty} \inf_{\mb J\in U \cap \mc B} \Phi_T(\mb J|\gamma) 
\end{equation}
Given $T>0$ there exist $\mb{J}^1\in U$ and $\mb W^1 \in \mc A_{T,\mb J^1}$
such that 
$$
\inf_{\mb J\in U} \Phi_T(\mb{J} |\gamma)  \ge 
\frac 1T I_{[0,T]}(\mb W^1|\gamma) - \frac 1T
$$
Let $\pi^1_t$, $t\in [0,T]$ the density path associated to $\mb W^1$,
$\delta>0$ and $T_0$ as in Lemma~\ref{t:cammino}.  
We now define, on the time interval $[0,T+T_0]$, the path 
$$
\mb W_t := \mb W^1_{t\wedge T} + \id_{[T,T+T_0]}(t) \mb W^2_{t-T}
$$ 
where $\mb W^2_t$, $t\in [0,T_0]$ is the path constructed in
Lemma~\ref{t:cammino} with $\gamma_1=\pi^1_T$ and $\gamma_2=\gamma$.
Let finally $\mb J= \mb W_{T+T_0} / (T+T_0)= T \mb J^1 /(T+T_0) +
\mb W^2_{T_0}/(T+T_0) $; note that, by construction,
$\mb J\in\mc B$. From Lemma~\ref{t:cammino} it now 
follows that, for $T$ large enough, $\mb J\in U$ and 
$$
\frac 1{T+T_0} I_{[0,T+T_0]}(\mb W|\gamma) \le 
\frac 1{T+T_0} I_{[0,T]}(\mb W^1|\gamma)
+ \frac 1{T+T_0} \big[ S_m(\gamma) +\delta \big]
\, .$$
By taking the limit $T\to\infty$, \eqref{toGlb} follows.
\end{proof}

\noindent\emph{Proofs of Proposition~\ref{t:Gcon} and Remark~\ref{t:gindip}.}
By Lemmata~\ref{t:d>d}, \ref{t:Gub}, and \ref{t:Glb} we have that, 
for each $\bar{\mb J} \in \mc M_d$ and any neighborhood $U$ 
of $\bar{\mb J}$ 
$$
\lim_{T\to \infty} \inf_{\mb J \in U} 
\Phi_T(\mb J |\gamma) 
= \inf_{\mb J \in U} {\tilde \Phi} (\mb J |\gamma)
$$
The $\Gamma$-convergence of the sequence $\Phi_T(\cdot |\gamma)$ to
the lower semi-continuous envelope of ${\tilde \Phi}(\cdot |\gamma)$ now
follows from the topological definition of  $\Gamma$-convergence, see
e.g.\ \cite[\S 1.4]{B}. 

We next prove  Remark~\ref{t:gindip}. 
Let $m\in (0,1)$; by using the path introduced in Lemma~\ref{t:cammino} it is
straightforward to check that, for each $\gamma_1,\gamma_2\in \mc F_m$
and $\mb J \in \mc B$ we have  
$\tilde \Phi(\mb J |\gamma_1)= \tilde \Phi(\mb J|\gamma_2)$ which
proves the first statement. 

Since $\Phi(\cdot|\gamma)$ is the lower semi-continuous envelope 
of $\tilde \Phi(\cdot|\gamma)$, it is enough to prove the convexity of
the latter. 
As $\mc B$ is a closed convex subset of $\mc M_d$, it is furthermore 
enough to show that for each $\mb J_1,\mb J_2 \in \mc B$, each $p\in
(0,1)$, and each $\gamma\in C^2(\bb T^d)$ bounded away from 0 and 1, 
we have 
\begin{equation}
\label{convex} 
\tilde \Phi(p\,\mb J_1 +(1-p) \,\mb J_2|\gamma)
\le p \, \tilde \Phi(\mb J_1|\gamma)
+ (1-p) \, \tilde \Phi(\mb J_2|\gamma)
\end{equation}
Given $\varepsilon>0$ we
can find $T>0$, $\mb W^1 \in \mc A_{pT,\mb J_1}$, and 
$\mb W^2 \in \mc A_{(1-p)T,\mb J_2}$ so that
\begin{eqnarray*}
\tilde \Phi(\mb J_1|\gamma) 
&\ge &
\frac {1}{pT} \: I_{[0,pT]} (\mb W^1|\gamma) -\epsilon
\\
\tilde \Phi(\mb J_2|\gamma) 
&\ge &
\frac 1{(1-p)T} \: I_{[0,(1-p)T]} (\mb W^2|\gamma) -\epsilon
\end{eqnarray*}
By the same arguments used in Lemma~\ref{t:subadd}, the path obtained
by gluing $\mb W^1$ with $\mb W^2$ is in the set 
$\mc A_{T, p \, \mb{J}_1+(1-p) \, \mb{J}_2 }$. 
The bound \eqref{convex} follows.
\qed

We conclude this section proving the exponential tightness needed to
complete the proof of Theorem~\ref{t:LDphi}.

\begin{lemma}
\label{gb1}
Fix a sequence $\eta^N\in \mc X_N$. There exists a sequence of compact
sets $\{K_\ell : \ell \ge 1\}$ of $\mc M_d$ such that
\begin{equation*}
\varlimsup_{T\to\infty} \varlimsup_{N\to\infty} 
\frac 1{T N^d} \log \bb P_{\eta^N}^N
\Big[ \frac 1T {\bf W}^N_T \in K_\ell^c \Big] \;\le\; - \ell \;. 
\end{equation*}
\end{lemma}

\begin{proof}
Fix a vector field $\mb H: \bb T^d \to \bb R^d$. We claim that
for every $A> 0$,
\begin{equation}
\label{gb2}
\bb P_{\eta^N}^N \Big[ \, \big\vert \<\mb W^N_T, \mb H\> \big\vert 
\ge A T  \Big] \;\le\; 2 \exp \big\{- TN^d [ A C(\mb H)^{-1} - 7] \big\}\;,
\end{equation}
where $C (\mb H) = \max_j \{ 1 \vee \Vert H_j \Vert_\infty^2 \vee
\Vert \partial_{u_j} H_j \Vert_\infty \}$.  Indeed, by Chebychev
exponential inequality,
\begin{equation*}
\bb P_{\eta^N}^N \Big[ \,  \<\mb W^N_T, \mb H\> 
\ge A T  \Big] \;\le\; e^{- \theta A T N^d}
\bb E_{\eta^N}^N \Big[ \mc M_T(\theta, \mb H)
e^{ f_N (T,\eta, \mb H) } \Big]
\end{equation*}
for every $\theta>0$.  Here, $\mc M_t(\theta, \mb H)$ is the mean one
exponential martingale defined just after \eqref{f12}, with a now time
independent vector field $\mb H$, and
\begin{eqnarray*}
f_N (T,\eta, \mb H) 
&=& N^2 \sum_{j=1}^d\sum_{x\in\bb T_N^d} \int_0^T ds\, 
\eta_s(x) [1-\eta_s(x+e_j)]  
\Big\{ e^{ \theta N^{-1} H_j (x/N)}-1 \Big\}   \\
&+& N^2 \sum_{j=1}^d\sum_{x\in\bb T_N^d} \int_0^T ds\, 
\eta_s(x+e_j) [1-\eta_s(x)] 
\Big\{ e^{ - \theta N^{-1} H_j (x/N)} - 1 \Big\} \;.
\end{eqnarray*}
A Taylor expansion shows that the absolute value of 
$f_N (T,\eta, \mb H)$ is bounded by 
$C (\mb H) T N^d \{ \theta + 2 \theta^2 e^{\theta C(\mb
  H)}\}$. To conclude the proof of the claim it remains to choose
$\theta = C(\mb H)^{-1}$, to remind that $\mc M_t(\theta, \mb H)$ has
mean one and to repeat the same argument with $-\mb H$ in place of
$\mb H$.

Recall the definition of the sequence $\{\mb G_k : k\ge 1\}$ defined
just after Proposition \ref{t:Gcon} and assume, without loss of
generality, that $C (\mb G_k) \le k$. For each $\ell \ge 1$, the set
$K_\ell$ of measures defined by 
\begin{equation*}
K_\ell \;=\; \bigcap_{k\ge 1} \{\mb J : \big\vert \<\mb J, \mb G_k\> 
\big\vert \le (k+7)^2 \ell \} 
\end{equation*}
is compact. On the other hand, by \eqref{gb2}, 
\begin{equation*}
\bb P_{\eta^N}^N \Big[\frac 1T {\bf W}^N_T \in K_\ell^c \Big] \;\le\;
4 e^{-TN^d \ell}
\end{equation*}
provided $N$ is sufficiently large. This proves the Lemma.
\end{proof}

\section{Dynamical phase transitions}
\label{app}

In this Section we analyze the variational problem \eqref{tilphi}
defining the functional $\tilde \Phi$. For the symmetric simple
exclusion process we prove, in Subsection~\ref{s:p=u}, 
that $\tilde \Phi=U$, where $U$ is defined in \eqref{dU}. 
Therefore no dynamical phase transition occurs in this case. 
In Subsection~\ref{s:dft} we consider a system with general transport
coefficients and show that, under suitable convexity assumptions, it
is possible to construct a traveling wave whose cost is, for $J$
large, strictly less than the constant profile. These convexity
hypotheses are satisfied for the KMP model \cite{bgl,kmp} and
therefore we prove that a dynamical phase transition takes place. 

\subsection{Symmetric simple exclusion process}
\label{s:p=u}

The following statement is essentially proven in \cite{bl6}; for the
reader's convenience we reproduce below its proof in a more formal
setting. Together with Theorem~\ref{t:LDphi} concludes the proof of 
Theorem~\ref{ldpnt}.

\begin{proposition}
\label{u=phi}
For each $m\in (0,1)$ the functional $U_m:\mc M_d \to [0,\infty]$ defined
in \eqref{dU} is lower semi-continuous. Moreover if $\gamma\in C^2(\bb
T^d)\cap \mc F_m$ is bounded away from 0 and 1 we have 
$U_m (\mb J)= \Phi (\mb J|\gamma)$ for any $\mb J\in \mc M_d$.
\end{proposition}

\begin{proof}
We first prove the lower semi-continuity of $U_m$.
Given $(\mb J,\rho) \in \mc M_d \times \mc F_m$ we define
$$
\mc U(\rho, \mb J) \; := \; \sup_{\mb F}  \mc V_{\mb F} (\rho,\mb J )
$$
where
$$
{\mc V}_{\mb F} (\rho,\mb J ) \; := \;   
\langle {\mb J} , \mb F \rangle  - \frac 12 
\langle \rho, \nabla \cdot \mb F \rangle 
- \frac 12 \langle  \mb F, \chi(\rho) \mb F \rangle
$$
and the supremum is carried over all smooth vector fields 
$\mb F \in C^1 (\bb T^d;\bb R^d)$. Note that, if $\mb J(du)= \mb j \, du$
for some $\mb j \in L_2(\bb T^d;\bb R^d)$ and $\rho \in C^2(\bb T^d)$
is bounded away from 0 and 1, we have
$$
\mc U(\rho, \mb J) = \frac 12 
\Big\langle \big[ {\mb j} + \frac 12 \nabla \rho \big], \frac 1{ \rho
    (1-\rho)} \big[ {\mb j} + \frac 12 \nabla \rho \big]
\Big\rangle
$$
Recalling definition \eqref{dU}, by the approximation arguments in
Subsection~\ref{sec4},  we have 
$U_m (\mb J) = \inf_{\rho \in \mc F_m} \mc U(\rho,\mb J)$. 

By the concavity of $\chi(\rho)$ we have that,
for each fixed $\mb F$, the functional 
$\mc V_{\mb F}(\cdot,\cdot) : \mc F \times \mc M_d \to \bb R$ is
convex and lower semi-continuous.  The lower semi-continuity of $\mc U$,
hence of $U_m$, follows now easily. We also note that the previous
argument shows that $\mc U$ is a convex functional on 
$\mc F\times \mc M_d$.

We next prove that for each $\gamma\in C^2(\bb T^d)\cap \mc F_m$
bounded away from 0 and 1 we have $U_m(\mb J) = \tilde \Phi (\mb
J|\gamma)$. From the definitions \eqref{dU}, \eqref{tilphi} and
Lemma~\ref{t:l2b}, we can assume that $\mb J(du) = \mb j \, du $ for
some $\mb j \in L_2(\bb T^d;\bb R^d)$ divergence free.

We first show that for each $T>0$, and each path 
$(\mb W, \pi)$ such that $\mb W\in \mc A_{T, \mb J}$ we have 
$$
\frac 1T I_{[0,T]} ( (\mb W, \pi) |\gamma ) \ge U_m(\mb J)
$$
Indeed, thanks to the approximation constructed in Section~\ref{sec4} we can
assume that $\pi$ is a smooth path bounded away from 0 and 1. 
For such a smooth path \eqref{f06} yields 
$$
\frac 1T I_{[0,T]} \big( (\mb W,\pi) | \gamma \big) = 
\frac 1T \int_0^T \! dt \: \mc U (\pi_t,\dot{\mb W}_t) 
\ge 
\mc U \Big( \frac 1T \int_0^T \!dt \: \pi_t \, , \, \mb J \Big)
\ge U_m (\mb J )
$$
where we used, in the second step, the joint convexity of the
functional $\mc U$ and Jensen inequality. 
In the last step we finally used that, by
conservation of mass $\pi_t \in \mc F_m$ for any $t\in [0,T]$. 

To show the converse inequality it is enough to construct, for
each $T$ large enough, an appropriate path.  
Given $\varepsilon >0$ there exists
$\rho\in C^2(\bb T^d)\cap \mc F_m$ bounded away from 0 and 1 such that 
$U_m(\mb j\, du ) \ge \mc U (\rho, \mb j) - \epsilon$. For $T>2 $ we 
construct the path  $(\mb W, \pi)$ such that
\begin{equation*}
\dot{\mb W}_t(du) = \left\{
\begin{array}{ccl}
\hat{\mb{w}} \, du &\textrm{if} & t \in [0,1) \\
\vphantom{\Big\{}
\frac{T}{T -2} \: \mb j \, du &\textrm{if} & t \in [1,T-1] \\
- \hat{\mb w} \, du &\textrm{if} & t\in (T-1,T]  
\end{array}
\right.
\end{equation*}
where $\hat{\mb w}$ solves $\nabla\cdot\hat{\mb w} = \gamma -\rho$.
It exists because $\gamma,\rho\in\mc F_m$, i.e.\ they have the same mass.
The density path $\pi$ is the corresponding solution of
\eqref{f09}, i.e.\ 
\begin{equation*}
\pi_t = \left\{
\begin{array}{ccl}
\gamma (1-t) + \rho t 
&\textrm{if} &  t\in [0,1) \\
\vphantom{\Big\{}
\rho   &\textrm{if} &  t \in [1, T -1] \\
\rho (T-t) + \gamma (T+1-t)
&\textrm{if} & t \in (T-1,T] 
\end{array}
\right.
\end{equation*}
It is straightforward to verify that $\mb W \in \mc A_{T, \mb J}$. 
Moreover 
$$
\lim_{T\to\infty} 
\frac 1T I_{[0,T]}\big( (\mb W, \pi) | \gamma \big) =
\mc U(\rho, \mb J)  \le U_m (\mb J) +\epsilon
$$
which concludes the proof.
\end{proof}

We conclude this Section by showing that in the one dimensional case
the functional $U_m$ is given by \eqref{dud=1}

\begin{lemma}
\label{s9}
Let $m\in (0,1)$, $d=1$ and $J(du)=j\, du$ for some $j\in \bb R$.  
Then
\begin{equation*}
U_{m}(J) \;=\; \frac 12 \frac{j^2}{\chi(m)}
\end{equation*}
\end{lemma}

\noindent\emph{Remark.} As it will be apparent from the proof, this 
Lemma holds whenever the real function $\rho\mapsto 1/\chi(\rho)$
is convex.  

\begin{proof}
Let $\rho \in \mc C^2(\bb T) \cap F_m$ be bounded away from 0
and 1. We have
\begin{equation*}
\int_{\bb T} \! du\, \frac{ [ j +  (1/2)\rho'(u)]^2 }
{\chi(\rho(u))} 
\;=\; \int_{\bb T} \!du \, \frac{j ^2}{\chi(\rho(u))} 
\;+\; \int_{\bb T} du \, \frac{[(1/2)\rho'(u)]^2}{\chi(\rho(u))} 
\end{equation*}
because the cross term vanishes upon integration.  
By Jensen inequality,
$$
\int_{\bb T} \! du \, \frac{j^2}{\chi(\rho(u))} 
\ge \frac{j^2} {\chi(m)} \;.
$$
On the other hand, by considering the constant profile $\rho(u) = m$,   
we trivially have
$
U_m(j \, du) \le  (1/2) \,  j^2 / \chi(m)
$. 
The lemma is therefore proven.
\end{proof}

\subsection{Other models}
\label{s:dft}

The general structure of the hydrodynamic equation obtained for the
scaling limit of the empirical density for stochastic lattice gases with a
weak external field $\mb E$ has the form, see \cite{kl,sp} 
\begin{eqnarray*}
&&  \partial_t \rho + \nabla \cdot \dot{\mb W} (\rho) = 0 \\
&& \dot{\mb W} (\rho) = - \frac 12 D(\rho) \nabla \rho + \chi(\rho) \mb E
\end{eqnarray*}
where $D(\rho)$ is the diffusion coefficient and $\chi(\rho)$ is the
mobility. 

In this general context, for smooth profiles, we let
$$
\mc U (\rho, {\mb J} )  = \frac 12 
\Big\langle \big[ {\mb J}  - \dot{\mb W} (\rho) \big], \frac 1{\chi(\rho)} 
\big[\mb J - \dot{\mb W}(\rho)\big] \Big\rangle
$$
The integrated empirical current is expected to satisfy, 
see \cite{bl6} for a heuristic derivation, 
a large deviation principle with rate function
$$
I_{[0,T]} ( \mb W ) = 
\int_0^T\!dt \: \mc U (\pi_t, \dot{\mb W}_t  )   
$$
in which $\pi$ is obtained from $ \dot{\mb W}_t $ by solving the
continuity equation $\partial_t \pi + \nabla \cdot \dot{\mb W} = 0$.

\medskip
We analyze the variational problem \eqref{tilphi} in this general
setting and  show that, under some assumptions on $D(\rho),\chi(\rho)$, 
a time dependent strategy is more convenient than taking a density path $\pi$
constant in time, so that $\Phi <U $. 
For simplicity, we here discuss only the one dimensional case $d=1$
and assume that there is no external field, $\mb E=0$; see \cite[\S~6.2]{bl6} 
for more details.

Given a mass $m$ and $v\in \bb R$, let 
$\Psi_v: \bb R \to \bb R_+$ be defined by
\begin{equation}
\label{PSI}
\Psi_v (J) \;=\; \inf_\rho \frac 12 \int_0^1 du\, \frac{\{J + 
v\, [\rho(u) - m] - w (\rho(u)) \}^2} {\chi(\rho(u))} \;,
\end{equation}
where  $ w (\rho) = \dot{W} (\rho) = -(1/2) D(\rho) \nabla\rho +\chi(\rho) E$  
and the infimum is carried over the profiles $\rho$ of mass $m$, i.e.\
over $\mc F_m$.
It will be convenient to write the term $D(\rho) \nabla\rho$ as 
$ \nabla d(\rho)$, i.e.\ $d(\rho)$ is an antiderivative of $D(\rho)$.

We claim that for each $v\in \bb R$
\begin{equation}
\label{g02}
\Phi \;\le\; \Psi_v
\end{equation}
Indeed, consider a profile $\rho_0$ in $\mc F_m$.
Let $T=v^{-1}$ and set $\rho(t,u) = \rho_0(u-vt)$, $w(t,u)= J + v
[\rho_0(u-tv) -m]$ in the time interval $[0,T]$. An elementary
computation shows that the continuity equation holds and that the time
average over the time interval $[0,T]$ of $w(\cdot, u)$ is equal $J$.
In particular,
\begin{equation*}
\Phi(J) \;\le\; \frac 1T \int_0^T dt \, \mc U (\rho (t), w (t))\; .
\end{equation*}

On the other hand, it is easy to show by periodicity that the right
hand side is equal to
\begin{equation*}
\frac 12 \int_0^1 du\, \frac{\{J + 
v[\rho_0(u) - m] - w (\rho_0)\}^2} {\chi(\rho_0(u))} \;.
\end{equation*}
Optimizing over the profile $\rho_0$, we conclude the proof of
\eqref{g02}.

\medskip
We next show that, if the real function $\rho \mapsto 1/\chi(\rho)$ is
convex and $\chi''(m)> 0$  for some $0<m<1$,  then  $\Phi(J) < U (J)$
for $|J|$ sufficiently large. 

\smallskip
To prove the previous statement, we first note that,
in view of \eqref{g02} and of Lemma \ref{s9}, it is enough to show
that there exists $\lambda\in\bb R$ such that
\begin{equation}
\label{PIPPO}
\limsup_{|J|\to\infty} \frac{\Psi_{\lambda J} (J)}{J^2} 
<  \frac{1}{2 \chi(m)} \;\cdot
\end{equation}

Fix a mass $m$, a current $J$ and take $v= \lambda J$. For $\rho\in \mc F_m$,
by expanding the square we get that
\begin{eqnarray}
\label{g04}
\!\!\!\!\!\!\!\!\!\!\!\! &&
\int_0^1 du\, \frac{\{J + 
\lambda J [\rho - m] + (1/2)\nabla d (\rho) \}^2}
{\chi(\rho)} \\
\!\!\!\!\!\!\!\!\!\!\!\! && \quad
=\; J^2 \int_0^1 du\, \frac{\{ 1 + \lambda [\rho -m]  \}^2}
{\chi(\rho)} \; +\;
\frac 14 \int_0^1 du\, \frac{[\nabla d(\rho)]^2}
{\chi(\rho)}\; \cdot 
\nonumber
\end{eqnarray}
because the cross term vanishes.  Expand the square on the first
integral. Let $F(r) = F_{\lambda,m}(r)$ be the smooth function defined by
\begin{equation*}
F(r) = \frac{\{ 1 + \lambda [r -m] \}^2}{\chi(r)} \;\cdot
\end{equation*}
An elementary computation shows that 
\begin{equation*}
F''(m) = \frac 1{\chi(m)^3} \Big\{ 2 \chi(m)^2 \lambda^2
- 4 \chi(m) \chi'(m) \lambda + 2 \chi'(m)^2 - \chi(m) \chi''(m) \Big\}\;.
\end{equation*}
Let $\lambda= \chi'(m)/\chi(m)$. For this choice $F''(m)<0$. In particular,
we can choose a non constant profile $\rho (u)$ in $\mc F_m$
close to $m$ such that $F''(\rho(u)) <0$ for every $u$. Hence, by
Jensen inequality, the coefficient of $J^2$ in \eqref{g04} is strictly
less than $\chi(m)^{-1}$.  This completes the proof of the claim.

\medskip
For the KMP model \cite{bgl,kmp} we have $D(\rho)=1$ and 
$\chi(\rho)= \rho^2$.  In particular $\chi$ and $1/\chi$ are convex
functions. Hence, by the above results, we have   
$\Phi(J) < U (J)$ for all sufficiently large currents $J$.

For the weakly asymmetric exclusion process with large external field
a similar phemenon occurs. More precisely, as shown in \cite{bd2},
there exists a traveling wave whose cost is strictly less than the
constant (in space and time) profile.  A numerical computation
\cite{bd2} suggests also that the minimizer of the variational problem
\eqref{tilphi} is indeed a traveling wave.

\bigskip\bigskip
\noindent{\bf Acknowledgments.} 
The authors acknowledge the support of PRIN MIUR 2004\-028108 
and 2004015228. 
A.D.S.\ was partially supported by Istituto Nazionale di 
Alta Matematica.
C.L.\ acknowledges the partial support of
the John S.\ Guggenheim Memorial Foundation, FAPERJ and CNPq.


\begin{thebibliography}{99}

\bibitem{BDGJL2} 
Bertini L., De Sole A., Gabrielli D., Jona--Lasinio G., Landim C., 
{\sl Macroscopic fluctuation theory for stationary non equilibrium
state.}
J.\ Statist.\ Phys.\ {\bf 107}, 635--675 (2002).

  
\bibitem{bdgjl3} 
Bertini L., De Sole A., Gabrielli D., Jona--Lasinio G., Landim C., 
{\sl Large deviations for the boundary driven simple exclusion
  process.}
Math.\ Phys.\ Anal.\ Geom.\ {\bf 6}, 231--267 (2003).
Math. Phys., Analysis and Geometry {\bf 6}, 231-267, (2003). 
  
\bibitem{BDGJL5} 
Bertini L., De Sole A., Gabrielli D., Jona--Lasinio G., Landim C., 
{\sl Current fluctuations in stochastic lattice gases.} 
Phys.\ Rev.\ Lett.\ {\bf 94}, 030601 (2005).


\bibitem{bl6} 
Bertini L., De Sole A., Gabrielli D., Jona--Lasinio G., Landim C., 
{\sl Non equilibrium current fluctuations in stochastic lattice gases.} 
Preprint 2005, arXiv: cond-mat/0506664.

  
\bibitem{bgl} 
Bertini L., Gabrielli D., Lebowitz J.L., 
{\sl Large deviations for a stochastic model of heat flow.}  
J.\ Statist.\ Phys.\ (2005). Available on line. 

\bibitem{bd} 
Bodineau T., Derrida B., 
{\sl Current fluctuations in non-equilibrium diffusive systems: an
  additivity principle.}  
Phys.\ Rev.\ Lett.\ {\bf 92}, 180601 (2004).

\bibitem{bd2} 
Bodineau T., Derrida B., 
{\sl Distribution of current in non-equilibrium diffusive systems and phase
transitions.} 
Preprint 2005, arXiv: cond-mat/0506540.


\bibitem{B} 
Braides A.,
{\sl $\Gamma$-convergence for beginners.} 
Oxford University Press, Oxford 2002.


\bibitem{dv} Donsker M.D., Varadhan S.R.S, 
{\sl Large deviations from a hydrodynamic scaling limit.}
Comm.\ Pure Appl.\ Math.\  {\bf 42}, 243--270 (1989).

\bibitem{kl} Kipnis C., Landim C., 
{\sl Scaling Limits of Interacting Particle Systems.}
Springer-Verlag, Berlin 1999.
  
\bibitem{kmp} Kipnis C., Marchioro C., Presutti E., 
{\sl Heat flow in an exactly solvable model.}  
J.\ Statist.\ Phys.\ {\bf 27}, 65--74 (1982).

\bibitem{kov} 
Kipnis C., Olla S., Varadhan S.R.S., 
{\sl Hydrodynamics and large deviations for a simple exclusion
  process.}
Comm.\ Pure Appl.\ Math.\ {\bf 42}, 115--137 (1989).


\bibitem{sp} 
Spohn H., 
{\sl Large scale dynamics of interacting particles.}
Springer-Verlag, Berlin 1991.



\end{thebibliography}
\end{document}